\pdfoutput=1
\documentclass[11pt,a4paper,twoside]{article}
\usepackage[T1]{fontenc}
\usepackage{lmodern,amsmath,amssymb,amsthm,mathtools,mathrsfs,microtype}
\usepackage[margin=24mm,headheight=12pt,headsep=18pt]{geometry}
\usepackage{fancyhdr,titlesec,enumitem,needspace,booktabs,tabularx}
\usepackage[hidelinks]{hyperref}
\titleformat{\section}{\centering\normalfont\scshape}{\thesection.}{.6em}{}
\titlespacing*{\section}{0pt}{18pt}{9pt}
\setlist[enumerate]{leftmargin=2em,itemsep=2pt,topsep=4pt}
\allowdisplaybreaks[1]\emergencystretch=1em
\AtBeginDocument{%
  \setlength{\abovedisplayskip}{10pt plus 2pt minus 4pt}%
  \setlength{\belowdisplayskip}{10pt plus 2pt minus 4pt}%
}
\newcommand{\OO}{\mathcal O}\newcommand{\II}{\mathcal I}
\newcommand{\QQ}{\mathcal Q}
\newcommand{\EE}{\mathcal E}\newcommand{\ZZ}{\mathcal Z}
\newcommand{\CC}{\mathcal C}\newcommand{\FF}{\mathcal F}
\newcommand{\VV}{\mathcal V}
\newcommand{\PP}{\mathbf P}
\newcommand{\C}{\mathbf C}
\newcommand{\Z}{\mathbf Z}\newcommand{\R}{\mathbf R}
\DeclareMathOperator{\Ext}{Ext}\DeclareMathOperator{\Hom}{Hom}
\DeclareMathOperator{\coker}{coker}
\DeclareMathOperator{\Pic}{Pic}\DeclareMathOperator{\Sing}{Sing}
\DeclareMathOperator{\Fitt}{Fitt}\DeclareMathOperator{\length}{length}
\DeclareMathOperator{\rk}{rk}
\DeclareMathOperator{\ch}{ch}
\DeclareMathOperator{\Hilb}{Hilb}
\DeclareMathOperator{\codim}{codim}\DeclareMathOperator{\pr}{pr}

\newcommand{\shExt}{\mathcal{E}\!xt}\newcommand{\shHom}{\mathcal H\!om}
\theoremstyle{plain}
\newtheorem{lemma}{Lemma}[section]
\newtheorem{theorem}[lemma]{Theorem}
\newtheorem{proposition}[lemma]{Proposition}
\newtheorem{corollary}[lemma]{Corollary}
\newtheorem{lettertheorem}{Theorem}

\theoremstyle{definition}
\newtheorem{definition}[lemma]{Definition}
\newtheorem{example}[lemma]{Example}

\theoremstyle{remark}
\newtheorem{remark}[lemma]{Remark}
\newcommand{\etale}{\textnormal{\'etale}}
\usepackage{etoolbox}
\AtBeginEnvironment{lettertheorem}{\Needspace{4\baselineskip}}
\AtBeginEnvironment{proposition}{\Needspace{5\baselineskip}}
\pretocmd{\section}{\Needspace{4\baselineskip}}{}{}
\pretocmd{\subsection}{\Needspace{4\baselineskip}}{}{}

\hypersetup{pdftitle={Families of rank-two reflexive sheaves on P1-bundles over K3 surfaces},pdfauthor={Qi Sun},pdfsubject={Stable reflexive sheaves, moduli components and Brill-Noether jump schemes}}
\begin{document}\thispagestyle{plain}
\begingroup
\renewcommand{\thefootnote}{\fnsymbol{footnote}}
\begin{center}
{\Large\bfseries Families of rank-two reflexive sheaves\\ on $\PP^1$-bundles over K3 surfaces}\par\vspace{7pt}
{\large Qi Sun\footnotemark[1]}\par
\end{center}
\footnotetext[1]{School of Mathematical Sciences, Shanghai Key Laboratory of PMMP,
East China Normal University, Math. Bldg, No. 500, Dongchuan Road, Shanghai,
200241, P. R. China, \href{mailto:52285500003@stu.ecnu.edu.cn}{52285500003@stu.ecnu.edu.cn}}
\endgroup
\begin{abstract}
We study deformations of stable rank-two reflexive sheaves on
$\PP^1$-bundles over K3 surfaces obtained by the Hartshorne--Serre
construction. The cokernel of a natural restriction map measures the
excess tangent dimension of the resulting extension family and gives a
criterion for its closure to be a generically smooth component of the
moduli space; enlarging the curve quotients via a relative Quot scheme
produces further generically smooth components, including nonsplit
examples. For split bundles $\OO_S\oplus\OO_S(eH)$ with
$\Pic(S)=\Z H$, $H^2=2h$, $e\ge2$ and $b\ge e+1$, we construct a
boundary divisor whose general point is stable and reflexive, with
unchanged $h^0$ but jumping tangent dimension. There the local
Brill--Noether ideal is nonzero and nilpotent, so the moduli space is
nonreduced although the number of sections does not jump. The divisor
also lies on a second component of dimension $(4b+2e+2)h$, formed by
extensions along torsion quotients, and the local dimension of the
moduli space at a general closed point of the divisor is exactly
$(4b+2e+2)h$.
\end{abstract}
{\footnotesize\noindent
\par}

\section{Introduction}
Let $S$ be a complex projective K3 surface, $\VV$ a rank-two vector
bundle on $S$, and $X=\PP_S(\VV)$ the associated $\PP^1$-bundle
(or threefold scroll over $S$).
For a rational polarization $L$, let $M_{X,L}(w)$ be the moduli
space of $L$-Gieseker-stable torsion-free sheaves on $X$ with fixed
numerical Chern character $w$. One approach to
its geometry is to study the Brill--Noether (BN) loci
\[
W^k_{X,L}(w)=\{[F]\in M_{X,L}(w):h^0(X,F)\ge k\},
\]
equipped with their determinantal scheme structures. Besides
nonemptiness and dimension, natural questions concern their
irreducible components, smoothness, and behavior under specialization.
These questions also involve the scheme structure of the Brill--Noether loci.

Theorem~\ref{thm:C} constructs stable reflexive sheaves $F_0$ for which
$h^0(X,F)$ is constant at all geometric points of an ambient neighborhood
of $[F_0]$, yet the local Brill--Noether ideal is nonzero and nilpotent,
detecting nonreducedness of the moduli space. For the split scroll
$X=\PP_S(\OO_S\oplus\OO_S(eH))$ with $\Pic(S)=\Z H$,
$H$ ample and base-point-free, $H^2=2h$ and $e\ge2$, this ideal has
$2eh-1$ generators with linearly independent images in the cotangent
space. These sheaves are dense in a boundary divisor shared by a generically
smooth component and a second component of larger dimension; we determine
the exact local dimension at a general closed point of this divisor.

For ruled surfaces, Costa and Mac\'ias Tarr\'io \cite{RuledSurfaces}
investigate nonemptiness and smoothness of Brill--Noether loci using
stability and the wall-and-chamber structure. On ruled threefolds over
$\PP^2$ \cite{RuledThreefolds}, they combine this approach with
Hartshorne--Serre extensions to describe components of moduli spaces
of stable rank-two bundles and to show that certain Brill--Noether
loci are nonempty. In this paper we consider
rank-two reflexive sheaves on $\PP^1$-bundles over K3 surfaces.
Our stability results in Section~\ref{sec:stability} follow a similar
strategy, but our main concern is different: we compute the tangent
correction to the Serre family, realize the missing directions by a
relative Quot construction, and study the scheme structure of the
Brill--Noether loci at the boundary of a generically smooth component. The Hartshorne--Serre
correspondence constructs such sheaves from curves and extension
sections \cite{Hartshorne,Arrondo}; their deformation theory and its
relation with Hilbert schemes are studied in
\cite{VermeireStable,VermeireModuli,Kleppe}. The sheaves in our
families have $c_3>0$, so they are not locally free; their deformation
theory must account for their isolated singularities.

Nonreduced parameter spaces already occur in Mumford's Hilbert-scheme
example \cite{Mumford}. Kleppe constructs nonreduced components of
moduli spaces of stable rank-two reflexive sheaves on $\PP^3$
\cite[Examples~2.4(b) and~2.5(a)--(b)]{Kleppe}; components and boundary loci for
rank-two torsion-free sheaves on $\PP^3$ are also studied in
\cite{JMT}. Gholampour and Kool \cite{GK} study stable reflexive
sheaves on $\PP^3$ using torus localization and compute generating
functions of Euler characteristics. Here nonreducedness is detected
by explicit local Brill--Noether ideals at points where the number of
sections is locally constant, complementing the deformation-theoretic
work of Kleppe and Vermeire.

Our aim is to determine when such an extension family fills a generically
smooth component and to study what happens at explicit reflexive
boundary points. Three questions guide the discussion: how can one
measure the tangent directions omitted by the section-induced construction;
can these directions be realized by enlarging its curve parameters;
and can Brill--Noether schemes detect nonreduced structure even when
the number of sections does not jump? Theorems~\ref{thm:A}--\ref{thm:C}
give answers for the families considered below.

Throughout, $H$ is an ample and base-point-free divisor on $S$, and
$\VV$ is globally generated. We use the quotient
convention on $\pi:X=\PP_S(\VV)\to S$, with
$\zeta=c_1(\OO_X(1))$, $A=\pi^*H$ and $H^2=2h=d$.
For the main component results assume $\Pic(S)=\Z H$ and
$\det\VV=\OO_S(eH)$. For $b>e/2$, write
\[
\xi_m=-\zeta+2mA,\qquad
\CC_b=\{L=\alpha\zeta+\beta A:\alpha,\beta>0,
\ \xi_bL^2<0<\xi_{b+1}L^2\}.
\]
We take polarizations to be rational classes in $\CC_b$.
For the fixed numerical Chern character considered below, the moduli
space is independent of this choice by Corollary~\ref{cor:chamber}.
The sheaves we construct are slope-stable, hence belong to
$M_{X,L}(w)$. Our determinantal convention for $W^k$ is recalled in
Section~\ref{sec:setup}; it does not require $H^2(F)=0$.

The quotient convention gives $\pi_*\OO_X(\zeta)=\VV$ and hence
$H^0(X,\OO_X(\zeta))\simeq H^0(S,\VV)$.
A nonzero $s\in H^0(S,\VV)$ defines a section of the line bundle
$\OO_X(\zeta)$; we write $D_s=(s=0)_X$ for its zero divisor on $X$,
rather than for the zero scheme of $s$ as a section of $\VV$ on $S$.
For a general pair $(s,C)$, with $C\in|H|$ smooth and $s|_C$ nowhere
zero, put $Z=D_s\cap\pi^{-1}(C)\simeq C$. The associated extensions are
\[
0\to\OO_X(bA)\to F\to\II_Z(\zeta-bA)\to0.
\]
For fixed $b$, we call the family of these middle terms, obtained by
varying the admissible pair $([s],C)$ and the nonzero extension class
up to scalar, the \emph{original Serre family}. Throughout this paper,
the term \emph{Serre family} denotes a family of middle terms of the
Hartshorne--Serre extensions specified in the corresponding construction.
Lemma~\ref{orig:6.1} constructs a flat
parameter family and its classifying map for the component results.
Set
\[
R_b=(4b+1)h,\qquad
D_b=(4b+2)h+h^0(\VV)-h^0(\VV(-H))-1,
\]
and let $\delta(s,C)$ be the dimension of the cokernel of
$H^0(S,\VV)\xrightarrow{s\wedge-}H^0(C,(\det\VV)|_C)$.

Concerning the first question, we prove that the restriction map
controls the difference between the parameter dimension and the
ambient tangent dimension.
\begin{lettertheorem}\label{thm:A}
Assume $\Pic(S)=\Z H$. For all sufficiently large $b$ and all rational
$L\in\CC_b$, the~image~$\FF_b$ of the original Serre family has
dimension $D_b$, and every member $F$ satisfies
\[
\dim\Ext^1_X(F,F)=D_b+\delta(s,C).
\]
Its reduced image closure is an irreducible component of
$M_{X,L}(w)_{\mathrm{red}}$ along which $M_{X,L}(w)$ is generically
smooth if and only if $\delta(s,C)=0$ for a general pair.
\end{lettertheorem}

Theorem~\ref{thm:A} is proved in Proposition~\ref{prop:5.7}, while
Theorem~\ref{thm:regular-component} covers every integer $b>e/2$ by
keeping track of the cohomology groups that vanish for $b\gg0$.

When $\delta>0$ for a general pair, the closure of the original
Serre family is not a component along which the ambient moduli
space is generically smooth. To address the second question,
we replace the curve quotients
$\VV|_C/\OO_Cs\simeq(\det\VV)|_C$ induced by global sections
with arbitrary line-bundle quotients on the curve having kernels of
degree zero. This gives the following result.
\begin{lettertheorem}\label{thm:B}
Assume that $\Pic(S)=\Z H$ and $e\ge2$, that $b$ is sufficiently
large, and that $L\in\CC_b$ is rational. Consider the family of all line-bundle quotients
$\VV|_C\twoheadrightarrow J$ on smooth $C\in|H|$ with kernel $K$ of degree zero.
The irreducible relative Quot component through a quotient induced
by an admissible pair $(s,C)$,
together with its Serre extensions, gives a generically smooth stable
moduli component of dimension
\[
D_b+\delta(s,C)=(4b+2e+1)h.
\]
Every member of the enlarged extension family has
$h^0(F)=b^2h+2+h^0(S,\VV(-(b+1)H))$.
\end{lettertheorem}

Proposition~\ref{orig:6.13} proves Theorem~\ref{thm:B}.
The correction $\delta(s,C)$ is the codimension in $\mathscr Q_0$ of
the closure of the quotient locus induced by global sections
(Remark~\ref{rem:delta-codim}).
The extra quotient directions need not all change $K$: restriction
of global sections to $C$ may itself fail to be surjective.
Example~\ref{orig:6.14} gives two nonsplit bundles on a double-plane K3 surface,
with $\delta=0$ and $\delta=1$, respectively.

Proposition~\ref{orig:6.12} computes sections and extension spaces for
quotients on smooth curves $C\in|H|$ with kernel of degree zero,
allowing torsion in $J$ as needed for the boundary construction and the
second component below.

We finally address the boundary question. In the split case, a
degeneration of the defining section produces stable reflexive
sheaves with unchanged $h^0$ but larger tangent spaces. More precisely,
take $\VV=\OO_S\oplus\OO_S(eH)$,
$e\ge2$ and $b\ge e+1$. Choose general
$u\in H^0(\OO_S(H))$, $v\in H^0(\OO_S(eH))$, and
$\tau\in H^0(\OO_S((2b+1)H))$ so that
$C=(u=0)$ is smooth and $\Gamma=C\cap(v=0)$ and
$\Xi=C\cap(\tau=0)$ are reduced and disjoint. Via the pushed-out
Koszul presentation of Lemma~\ref{orig:6.3}, the degenerate
section $s=(0,v)$ defines
\[
0\to\OO_X(-(b+1)A)\xrightarrow{(\tau,u,s)}
\OO_X(bA)\oplus\OO_X(-bA)\oplus\OO_X(\zeta-(b+1)A)
\to F_0\to0.
\]
Put $\ell=2eh$. For this split bundle and $b\ge e+1$, one has
$k_b=b^2h+2$ and $D_b=(4b+2e+1)h$.
We denote by $\mathcal M_b$ the reduced closure in
$M_{X,L}(w)$ of the original Serre family for this split bundle,
where $w=\ch(F_0)$ and $s$ has a nonzero constant component.

The next theorem computes these jumps and shows that the local
Brill--Noether ideal at $[F_0]$ is nonzero but contained in the nilradical.
\begin{lettertheorem}\label{thm:C}
Assume $\Pic(S)=\Z H$, $e\ge2$ and $b\ge e+1$, and let
$L\in\CC_b$ be rational.
\begin{enumerate}[label=\textup{(\roman*)}]
\item The sheaf $F_0$ above is a stable reflexive boundary point of the
generically smooth Serre component $\mathcal M_b$.
Its reduced singularity scheme has length $(4b+2)h$, and
\[
\begin{gathered}
h^0(F_0)=k_b,\qquad h^1(F_0)=\ell,\qquad h^3(F_0)=0,\\
h^2(F_0)=4+(e^2-2eb+2b^2+2b+1)h,\\
\dim\Ext^1_X(F_0,F_0)=D_b+(2e+1)h.
\end{gathered}
\]
\item In the full stable moduli space there is a neighborhood $\Omega$ of
$[F_0]$ with $W^{k_b+1}\cap\Omega=\varnothing$.
At $[F_0]$, the local ideals of $W^{k_b}$ and
$\mathcal B_b=\{[F]:h^0(F(-bA))\ge1\}$ coincide and lie in
the nilradical. Their common ideal has $\ell-1$ generators whose
images in the cotangent space are linearly independent.
In particular, the ambient moduli scheme is nonreduced at $[F_0]$.
\item As $(u,v,\tau)$ vary, the closure of these boundary points is an
irreducible divisor $\Delta_b\subset\mathcal M_b$
of dimension $D_b-1$, and \textup{(i)}--\textup{(ii)} hold at its general point.
\item There is an irreducible component $\mathcal N_b$ of
$M_{\mathrm{red}}$, distinct from $\mathcal M_b$, of
dimension $D_b+h$, such that
\[
\Delta_b\subset \mathcal M_b\cap \mathcal N_b ,
\]
and $\Delta_b$ is an irreducible component of this intersection. A general member of $\mathcal N_b$ is a stable reflexive extension
$0\to\mathcal O_X(bA)\to F\to\mathcal I_Z(\zeta-bA)\to0$ with
$Z=\PP_C(\mathcal O_C\oplus\mathcal O_{\Gamma'}(eH))$ for a smooth
$C\in|H|$ and a reduced divisor $\Gamma'\subset C$ of degree $\ell$.
At a general closed point $[F_0]\in\Delta_b$,
\[
\dim\mathcal O_{M,[F_0]}=D_b+h .
\]
\end{enumerate}
\end{lettertheorem}

Theorem~\ref{thm:C} follows from Propositions~\ref{orig:9.6},
\ref{cor:k3-boundary-nilpotents}, \ref{prop:torsion-sheaves}
and~\ref{prop:torsion-family} and Corollaries~\ref{cor:boundary-divisor}
and~\ref{cor:second-component}. Part (ii) detects local nonreducedness at $[F_0]$
through the BN ideals, without identifying the full nilradical. For the
split bundle, $\delta=0$, so the original family already fills the
component of Theorem~\ref{thm:B}. The phenomenon in Theorem~\ref{thm:C} arises upon
specialization to its boundary, where the Serre component meets a second
component $\mathcal N_b$. This component comes from the same
Hartshorne--Serre construction, but with torsion quotients: the kernel
$K$ of $\mathcal V|_C\twoheadrightarrow J$ is an arbitrary degree-zero
line bundle of the form $\mathcal O_C(eH-\Gamma')$ rather than
$\mathcal O_C$. The local BN ideal detects tangent directions that cannot
be induced by a reduced parameter space; part (iv) shows that the reduced
local dimension at $[F_0]$ is $D_b+h$, one less than the dimension
$D_b+h+1$ of the BN tangent space $\ker\rho$, where $\rho$ is the
section-lifting map of Proposition~\ref{orig:9.6}. We do not know whether
$\mathcal M_b$ and $\mathcal N_b$ are the only components
through a general point of $\Delta_b$, nor whether $M$ is generically
reduced along $\mathcal N_b$.

The paper is organized as follows. Section~\ref{sec:setup} fixes
the notation and recalls the cohomological and determinantal
conventions. Section~\ref{sec:construction} constructs the Serre
extensions and computes their singularities and the
cohomology of curve quotients. Section~\ref{sec:stability} establishes
an explicit stability interval for arbitrary Picard rank
and excludes additional numerical walls when $\Pic(S)=\Z H$.
Section~\ref{sec:components} develops the component and tangent-space
criteria used in Theorem~\ref{thm:A}.
Section~\ref{sec:quot} enlarges the quotient family, proves
Theorem~\ref{thm:B}, and gives nonsplit examples.
Section~\ref{sec:boundary} proves Theorem~\ref{thm:C} by constructing
the boundary sheaves $F_0$ and the divisor $\Delta_b$, analyzing the local
Brill--Noether ideals, and constructing a second component
$\mathcal N_b$ through $\Delta_b$ in Subsection~\ref{subsec:second-component}
to determine the exact local dimension.

\section{Notation and preliminaries}\label{sec:setup}

Throughout, $S$ is a smooth complex projective K3 surface, $H$ is ample
and base-point-free, and $\VV$ is globally generated of rank two. Put
$X=\PP_S(\VV)$ in the quotient convention, $A=\pi^*H$ and
$\zeta=c_1(\OO_X(1))$. The same quotient convention applies to
$\PP_C(J)$. For spaces of sections and bundles of extension classes,
$\PP(W)$ and $\PP_B(\EE)$ instead parametrize one-dimensional
subspaces, with tautological subbundle $\OO(-1)$.
Write $\gamma=\deg c_2(\VV)$ and $d=H^2=2h$.
We use $K_S\simeq\OO_S$, $H^1(S,\OO_S)=0$ and
$\chi(\OO_S)=2$ throughout.
For a coherent sheaf $E$ on a scheme $Y$, we abbreviate
$H^i(Y,E)$ and $h^i(Y,E)$ to $H^i(E)$ and $h^i(E)$ when the
base is clear. Numerical equivalence is denoted by $\equiv$.
The additional Picard-rank hypothesis
\[
\Pic(S)=\Z H \tag{N}\label{eq:N}
\]
is imposed only when explicitly stated; under~\eqref{eq:N}, we
write $\det\VV\simeq\OO_S(eH)$ with $e\ge0$.

Global generation and rank-two duality give a subbundle inclusion
\[
\VV\simeq\VV^\vee\otimes\det\VV
\hookrightarrow(\det\VV)^{\oplus h^0(S,\VV)}.
\]

We use the standard notions of torsion-free, reflexive and saturated
sheaves, with $F^\vee=\shHom(F,\OO_X)$. Rank-one reflexive sheaves
are line bundles; reflexive sheaves on a smooth threefold are locally
free off a finite set. We use the torsion-free $S_2$ criterion
\cite[\S1]{Hartshorne}. For our rank-two sheaves,
$\Sing(F)=V(\Fitt_2(F))$ records the scheme structure of the
non-locally-free locus.

The slope is $\mu_L(F)=c_1(F)\cdot L^2/\rk F$. Recall the
implications slope-stable $\Rightarrow$ Gieseker-stable $\Rightarrow$
slope-semistable. Gieseker stability is taken with respect to an integral
multiple of a rational polarization, while numerical walls are
considered in the cone of real ample classes. Let $M_{X,L}(w)$ denote the
moduli space of Gieseker-stable torsion-free sheaves of fixed numerical
Chern character $w$; superscripts $\mathrm{ref}$ and $\mathrm{vb}$
denote its reflexive and locally free open loci. All families are assumed flat over their parameter schemes.
The moduli stack of stable sheaves is a $\mathbf G_m$-gerbe over
$M_{X,L}(w)$ and admits universal families on \etale{} charts
\cite{HL,Lieblich}.
Scalar automorphisms act trivially on $\Ext^1(F,F)$, so
$T_{[F]}M_{X,L}(w)=\Ext^1_X(F,F)$.

The Chow relation and canonical bundle formula are
\[
\zeta^2-\zeta\pi^*c_1(\VV)+\pi^*c_2(\VV)=0,\qquad
K_X=-2\zeta+\pi^*c_1(\VV).\tag{2.1}\label{eq:orig:2.1}
\]
The second follows from the relative Euler sequence
$0\to\Omega^1_{X/S}(\zeta)\to\pi^*\VV\to\OO_X(\zeta)\to0$.
\begin{lemma}\label{orig:2.7}
One has $A^3=0$, $\zeta A^2=d$, $\zeta^2A=c_1(\VV)H$, and
$\zeta^3=c_1(\VV)^2-\gamma$ as degrees. Under \eqref{eq:N} these become
\[
A^3=0,\quad\zeta A^2=d,\quad\zeta^2A=ed,\quad\zeta^3=e^2d-\gamma,
\quad N^1(X)_\R=\R\zeta\oplus\R A.
\]
\end{lemma}
\begin{proof}
Degrees of zero-cycles are preserved by $\pi_*$, and the projection formula gives
$\pi_*(\pi^*\beta\cdot\alpha)=\beta\cdot\pi_*\alpha$.
The projective-bundle formula \cite[\S3]{Fulton} gives
$\pi_*(\zeta\cap[X])=[S]$, and $\pi_*[X]=0\in\operatorname{CH}_3(S)$
since $\dim S<\dim X$.
Applying $\pi_*$ to \eqref{eq:orig:2.1} and to its product with $\zeta$ gives
\[
\pi_*(\zeta^2)=c_1(\VV),\qquad
\pi_*(\zeta^3)=c_1(\VV)^2-c_2(\VV).
\]
Hence $A^3=\pi^*(H^3)=0$, $\zeta A^2=H^2=d$,
$\zeta^2A=c_1(\VV)\cdot H$, and $\zeta^3=c_1(\VV)^2-\gamma$ as degrees.
Under \eqref{eq:N}, we have $c_1(\VV)=eH$ and
$\Pic(X)=\Z[\OO_X(1)]\oplus\pi^*\Pic(S)=\Z\zeta\oplus\Z A$
by \cite[II, Exercise~7.9]{HartshorneAG}.
The classes $\zeta$ and $A$ are numerically independent: intersecting a numerical
relation $a\zeta+bA\equiv0$ first with $A^2$ and then with $\zeta A$ gives
$ad=0$ and $bd=0$, giving $N^1(X)_\R=\R\zeta\oplus\R A$.
Under \eqref{eq:N}, $\Theta=\pi^*c_2(\VV)-(\gamma/d)A^2$ pairs trivially
with both generators $A,\zeta$ of $N^1(X)_\R$, hence $\Theta\equiv0$ and
$\zeta^2\equiv e\zeta A-(\gamma/d)A^2$.
\end{proof}
Global generation makes $\zeta$ nef and relatively ample. A generating
space gives a closed immersion $X\hookrightarrow S\times\PP^m$,
so restriction of positive combinations of the two ample factors
shows that
\[
\mathcal A_+=\{\alpha\zeta+\beta A:\alpha,\beta>0\}
\subset\operatorname{Amp}(X)_\R.\tag{2.2}\label{eq:orig:2.2}
\]
We work throughout with this explicit ample subcone.

The cohomology of line bundles on $X$ reduces to the base as follows.
\begin{lemma}\label{orig:2.8}
For $m\geq0$ and every integer $t$,
\[
H^i(X,\OO_X(m\zeta+tA))
\simeq H^i(S,\operatorname{Sym}^m\VV\otimes\OO_S(tH)).
\]
For line bundles of relative degree $-1$, all the groups
$H^i(X,\OO_X(-\zeta+tA))$ vanish.
\end{lemma}
\begin{proof}
In the quotient convention, relative cohomology of projective bundles
\cite[II, Proposition 7.11; III, Exercise 8.4]{HartshorneAG} gives
$\pi_*\OO_X(m\zeta)\simeq\operatorname{Sym}^m\VV$ and
$R^q\pi_*\OO_X(m\zeta)=0$ for $m\geq0$, $q>0$.
In relative degree $-1$, $R^q\pi_*\OO_X(-\zeta)=0$ for every
$q\geq0$. The projection formula and the Leray spectral sequence
give the asserted cohomology identities for every integer $t$.
\end{proof}
The cohomology of $\OO_S(nH)$ has the following explicit form.
For every integer $n\ne0$,
\[
H^1(S,\OO_S(nH))=0,\qquad
h^0(S,\OO_S(nH))=
\begin{cases}
2+n^2h,&n>0,\\
0,&n<0.
\end{cases}
\]
Indeed, for $n>0$, Kodaira vanishing gives
$H^i(S,\OO_S(nH))=0$ for $i=1,2$, and Riemann--Roch gives
$h^0(S,\OO_S(nH))=2+n^2h$.
For $n<0$, Serre duality identifies $H^1(S,\OO_S(nH))$ with
$H^1(S,\OO_S(-nH))^\vee$, while $nH\cdot H<0$ excludes a
nonzero global section of $\OO_S(nH)$.
For $n=0$, $(h^0,h^1,h^2)=(1,0,1)$. Moreover,
$\chi(\OO_S(nH))=2+n^2h$ for all $n$, and
$H^2(S,\OO_S(nH))=0$ for $n>0$.
We write $\nu_n=h^0(S,\OO_S(nH))$, so that $\nu_n=2+n^2h$ for
$n>0$, $\nu_0=1$ and $\nu_n=0$ for $n<0$.
If $C\in|H|$ is smooth, the ideal sequence and
$H^1(S,\OO_S(-H))=0$ show that $C$ is connected. Adjunction gives
\[
K_C\simeq\OO_C(H),\qquad g(C)=h+1.
\]

\Needspace{5\baselineskip}
For section counts in families, we use the following determinantal
scheme structure on the Brill--Noether loci, defined by the full
cohomology complex \cite{BudurWang}. For the classical construction
for vector bundles under higher-cohomology vanishing, see
\cite[Theorem~2.3]{CostaMiroRoig}.
\begin{definition}\label{orig:2.9}\label{def:bn}
For $k\ge0$, the \emph{Brill--Noether scheme} $W^k_{X,L}(w)$ is
the determinantal subscheme of $M_{X,L}(w)$ with geometric points
$\{[F]:h^0(X,F)\ge k\}$.
On an \etale{} universal chart $T$, with flat family $\mathscr F$
and $p:X\times T\to T$, represent the full cohomology complex locally by
\[
Rp_*\mathscr F\simeq[K^0\xrightarrow{d^0}K^1\to K^2\to K^3],
\]
where the $K^i$ are finite-rank locally free sheaves in degrees
$0,\ldots,3$. Derived base change gives
$H^0(X_t,\mathscr F_t)=\ker(d^0\otimes\kappa(t))$
\cite[Theorem 2.5.4]{TT}, \cite[Exp.~IV, Proposition 3.1.0]{SGA6}.
For $r=\rk K^0$, the ideal of
$W^k_{X,L}(w)\times_{M_{X,L}(w)}T$ is
\[
I_{r-k+1}(d^0),
\]
where $I_j$ denotes the ideal of $j$-minors, with $I_j=\OO_T$
for $j\le0$ and $I_j=0$ when $j$ exceeds the matrix dimensions.
\end{definition}
These ideals are independent of the chosen representative, unchanged
by tensoring $\mathscr F$ with a line bundle from $T$, and compatible
with arbitrary base change; hence they descend to $M_{X,L}(w)$
without requiring any higher cohomology vanishing
\cite[Definition--Proposition 2.2 and Corollaries 2.3--2.5]{BudurWang}.

A saturated rank-one subsheaf can change stability only where its slope equals that of the rank-two sheaf. We encode this equality by a numerical divisor class and record the effectivity condition imposed by the second Chern class.

\begin{definition}\label{orig:2.11}\label{def:walls}
Fix numerical Chern classes $c_1,c_2$ of rank-two sheaves on $X$.
For an integral numerical divisor class
$0\ne\eta\equiv-c_1\pmod2$, put
\[
\Delta_\eta=c_2+\frac{\eta^2-c_1^2}{4},\qquad
W^\eta=\{L\in\operatorname{Amp}(X)_\R:\eta\cdot L^2=0\}.
\]
A nonempty $W^\eta$ is a \emph{numerical rank-one wall}. It is
\emph{numerically effective} if $\Delta_\eta$ can be represented
by an effective (possibly zero) codimension-two cycle.
\end{definition}
For a rank-two torsion-free sheaf $F$, write a saturated rank-one
subsheaf as $M_1=\II_T\otimes L_1$ and its torsion-free quotient as
$\II_Y\otimes L_2$, where $L_1,L_2$ are line bundles and
$\codim T,\codim Y\ge2$. Then $c_1(L_2)=c_1(F)-c_1(L_1)$, and the exact sequence
\[
0\to\II_T\otimes L_1\to F\to\II_Y\otimes L_2\to0
\]
gives $c_2(F)=c_1(L_1)\bigl(c_1(F)-c_1(L_1)\bigr)+[T]_2+[Y]_2$.
Hence, for $\eta=2c_1(L_1)-c_1(F)$,
\[
\Delta_\eta=c_2(F)-c_1(L_1)\bigl(c_1(F)-c_1(L_1)\bigr)=[T]_2+[Y]_2.
\]
If $F$ is reflexive, its saturated rank-one subsheaf $M_1$ is a line
bundle, and the term $[T]_2$ vanishes.
Numerical effectivity is necessary for a saturated destabilizing
sequence, but does not imply an actual stability change.

For these fixed Chern classes, a \emph{chamber} is a connected
component of the complement in $\mathcal A_+$ of the union of
the numerically effective rank-one walls.

The following descent lemma will be used for the universal family in
Section~\ref{subsec:local}.

\Needspace{4\baselineskip}
\begin{lemma}\label{lem:pullback-descent}
Let $\pi:X\to S$ be a smooth projective morphism of smooth projective
varieties over $\C$ such that the natural map
$\OO_S\to R\pi_*\OO_X$ is an isomorphism. Let $U$ be a scheme of
finite type over $\C$, $t_0\in U$ a closed point, and $\QQ$ a coherent
$U$-flat sheaf on $X\times U$ with $\QQ_{t_0}\simeq\pi^*\mathcal A_0$ for a
coherent sheaf $\mathcal A_0$ on $S$. Set $\pi_U=\pi\times\mathrm{id}_U$.
After shrinking $U$ around $t_0$, the sheaf
$\mathcal A:=\pi_{U*}\QQ$ is coherent and $U$-flat, and the adjunction map
\[
\pi_U^*\mathcal A\longrightarrow\QQ
\]
is an isomorphism. The formation of $\mathcal A$ commutes with every base change
$T\to U$: writing $\mathcal A_T$ and $\QQ_T$ for the pullbacks, the
canonical map $\mathcal A_T\to\pi_{T*}\QQ_T$ is an isomorphism.
\end{lemma}
\begin{proof}
This is a derived version of the seesaw principle
(cf.~\cite[\S5]{MumfordAV}).

The coherent sheaf $\QQ$ is perfect relative to $U$ by its
$U$-flatness. Since the projection $X\times U\to U$ is smooth,
$\QQ$ is also perfect on $X\times U$
\cite[Section~2.1]{Lieblich}.
Since $\pi_U$ is smooth and projective,
$K:=R\pi_{U*}\QQ$ is perfect and commutes with derived base change
\cite[Theorems~2.5.4 and~2.5.6]{TT}.
Let $\theta:\mathbf L\pi_U^*K\to\QQ$ be the counit.

By derived base change, its restriction to $X\times\{t_0\}$ is
\[
\mathbf L\pi^*R\pi_*\mathbf L\pi^*\mathcal A_0
\longrightarrow\mathbf L\pi^*\mathcal A_0.
\]
The projection formula identifies the unit
$\mathcal A_0\to R\pi_*\mathbf L\pi^*\mathcal A_0$ with
$\mathcal A_0\to \mathcal A_0\otimes^{\mathbf L}R\pi_*\OO_X$, which is an
isomorphism by hypothesis. The adjunction identity makes
the restricted counit an isomorphism. The perfect complex
$\operatorname{Cone}(\theta)$ has zero derived restriction to the
central fiber. By derived Nakayama \cite[Lemma~2.1.3]{Lieblich},
its closed support is disjoint from that fiber. Its image in $U$ is closed
because $X$ is proper. Removing this image makes $\theta$ an isomorphism.

Let $t\in U$ be a closed point and set
$K_t:=K\otimes_{\OO_U}^{\mathbf L}\kappa(t)$.
Base change gives $\QQ_t\simeq\mathbf L\pi^*K_t=\pi^*K_t$,
since $\pi$ is flat. The equality $\pi_*\OO_X=\OO_S$ implies
that $\pi$ is surjective, hence faithfully flat.
\[
\pi^*\mathcal H^i(K_t)\simeq\mathcal H^i(\QQ_t)=0
\quad(i\ne0)
\]
forces $K_t$ to be a sheaf. Applying the fiberwise sheaf and
flatness criterion \cite[Lemma~4.3, corrected version]{Bridgeland}
(see also \cite[Lemma~3.31]{HuybrechtsFM}) to the flat projection
$S\times U\to U$, we obtain $K\simeq\mathcal A[0]$ for a
coherent $U$-flat sheaf $\mathcal A$.
We have $\mathcal A=\pi_{U*}\QQ$ and
$\pi_U^*\mathcal A\simeq\QQ$.

Finally, for any $u:T\to U$, the $U$-flatness of $\mathcal A$
and $\QQ$ and derived base change give
\[
\mathcal A_T\simeq\mathbf L(\mathrm{id}_S\times u)^*K
\simeq R\pi_{T*}\QQ_T.
\]
The derived pullbacks here are ordinary pullbacks because both sheaves
are flat over $U$; the isomorphism is the natural base-change morphism.
Taking $\mathcal H^0$ gives the canonical ordinary base-change
isomorphism claimed in the statement.
\end{proof}

We will also use the following criterion to recognize a moduli component
from a family whose dimension equals the ambient tangent dimension.

\begin{lemma}\label{lem:regularity}
Let $M$ be a scheme of finite type over $\C$, let $Y\subset M$ be an
irreducible closed subset, and let $y\in Y$ be a closed point such that
$\dim_\C T_yM\le\dim Y$. Then $\OO_{M,y}$ is regular, $Y$ is the unique
irreducible component of $M$ through $y$, and $M$ is smooth along a
dense open subset of $Y$.
\end{lemma}
\begin{proof}
Since $y$ is closed and $Y$ is irreducible of finite type over $\C$,
$\dim\OO_{Y,y}=\dim Y$. Hence
$\dim Y\le\dim\OO_{M,y}\le\dim_\C T_yM\le\dim Y$.
Equality of dimension and embedding dimension makes $\OO_{M,y}$
regular, hence a domain. The unique component through $y$ contains
$Y$ and has the same dimension, so it equals $Y$.
Over $\C$, regularity at $y$ implies smoothness, which is open.
\end{proof}

\section{Serre extensions and section counts}\label{sec:construction}
Throughout this section, we fix an integer $b\ge1$. We first construct a smooth complete-intersection
curve $Z\subset X$ isomorphic to a member of $|H|$ and compute its
extension space. We then determine the local structure and Chern
classes of the resulting sheaves. Theorem~\ref{thm:serre-ext}
collects these results; their stability is proved in
Section~\ref{sec:stability}.

We consider extensions
\[
0\to\OO_X(bA)\to F_\sigma\to\II_Z(\zeta-bA)\to0.\tag{3.1}\label{eq:orig:4.1}
\]

We first choose the complete-intersection curve so that its projection
to the base is an isomorphism.
By global generation and generic smoothness, a general section
$s\in H^0(S,\VV)$ has a finite reduced (possibly empty) zero scheme
$(s)_0$ of length $\gamma$ \cite[III, Corollary 10.7]{HartshorneAG}.
By Bertini's theorem, a general member $C=(u=0)\in|H|$, defined
by $u\in H^0(S,\OO_S(H))$, is smooth and disjoint from $(s)_0$
\cite[III, Corollary 10.9]{HartshorneAG}.
The vanishing of $H^1(S,\OO_S(-H))$ recalled in Section~\ref{sec:setup}
shows that $C$ is connected.

Under $H^0(S,\VV)\simeq H^0(X,\OO_X(\zeta))$, the section $s$
induces a section of the tautological line bundle on $X$. We denote
its zero divisor by $D_s=(s=0)_X$, to distinguish it from $(s)_0$,
and put $Z=D_s\cap\pi^{-1}(C)$. Since $C\cap (s)_0=\varnothing$,
the induced section on $\pi^{-1}(C)$ is a nonzero linear form on
every fiber.
The quotient $\VV|_C/\OO_Cs\simeq(\det\VV)|_C$ identifies
\[
Z=\PP_C((\det\VV)|_C)\simeq C,\qquad
\OO_Z(\zeta)\simeq(\det\VV)|_C,\quad\OO_Z(A)\simeq\OO_C(H). \tag{3.2}\label{eq:orig:4.2}
\]
Since $s$ cuts out a Cartier divisor on $\pi^{-1}(C)$, the sections
$\pi^*u$ and $s$ form a regular sequence on $X$. Thus $[Z]=\zeta A$,
and the Koszul resolution is
\[
0\to\OO_X(-\zeta-A)\to\OO_X(-\zeta)\oplus\OO_X(-A)
\to\II_Z\to0. \tag{3.3}\label{eq:orig:4.3}
\]

With the curve fixed, we next identify the full extension space with a
space of sections on the curve. This will allow us to describe the
singularities of an extension in terms of the zeros of the corresponding
section.

\begin{lemma}\label{orig:4.2}
For every $b\geq1$ there is a natural isomorphism
\[
\Ext^1_X(\II_Z(\zeta-bA),\OO_X(bA))
\simeq H^0(Z,\OO_Z((2b+1)A)).
\]
Under \eqref{eq:orig:4.2}, this is $H^0(C,K_C^{\otimes(2b+1)})$.
Moreover,
\[
\begin{aligned}
\dim\Ext^1_X(\II_Z(\zeta-bA),\OO_X(bA))&=R_b=(4b+1)h,\\
\Ext^2_X(\II_Z(\zeta-bA),\OO_X(bA))&=0.
\end{aligned}
\]
\end{lemma}
\begin{proof}
On $X\setminus Z$, the ideal $\II_Z$ equals $\OO_X$, so a morphism
$\II_Z(\zeta-bA)\to \OO_X(bA)$ is a section of $\OO_X(2bA-\zeta)$ there.
Since $X$ is normal and $Z$ has codimension two, this section extends
uniquely across $Z$. Conversely, multiplication by a section of
$\OO_X(2bA-\zeta)$ gives a morphism $\II_Z(\zeta-bA)\to \OO_X(bA)$.
We obtain $\shHom(\II_Z(\zeta-bA),\OO_X(bA))\simeq \OO_X(2bA-\zeta)$. Locally the resolution
$0\to R\to R^2\to(f,g)\to0$ gives
$\Ext^1_R((f,g),R)=R/(f,g)$ and vanishing in degrees at least two.
Local duality for the codimension-two embedding $Z\hookrightarrow X$
in the Hartshorne--Serre correspondence \cite[\S1]{Arrondo} gives
\[
\shExt^1_X(\II_Z(\zeta-bA),\OO_X(bA))\simeq\det N_{Z/X}\otimes \OO_Z(2bA-\zeta).
\]
Since $\det N_{Z/X}=\OO_Z(\zeta+A)$, this is
$\OO_Z((2b+1)A)$. All cohomology of
$\OO_X(2bA-\zeta)$ vanishes by Lemma~\ref{orig:2.8}.
The local-to-global Ext spectral sequence gives
\[
\begin{gathered}
\Ext^i_X(\II_Z(\zeta-bA),\OO_X(bA))\simeq
H^{i-1}(Z,\OO_Z((2b+1)A))\quad(i\geq1),\\
\Hom_X(\II_Z(\zeta-bA),\OO_X(bA))=0.
\end{gathered}
\]
Indeed, the $E_2$-page is concentrated in the row $q=1$,
corresponding to the sheaf $\shExt^1_X(\II_Z(\zeta-bA),\OO_X(bA))$ supported on $Z$.
This identifies the full extension space.

Since $K_C\simeq\OO_C(H)$, Serre duality gives
\[
H^1(C,\OO_C((2b+1)H))^\vee=H^0(C,\OO_C(-2bH))=0
\]
for $b\ge1$. Curve Riemann--Roch then gives
$R_b=(2b+1)\cdot2h+1-(h+1)=(4b+1)h$ and the claimed
degree-two Ext vanishing.
\end{proof}

We next examine the local extension class to determine the singularities
and establish reflexivity.

\begin{lemma}\label{orig:4.3}
Let $(R,\mathfrak m)$ be a regular local ring of dimension three and let
$I=(f,g)$ be generated by a regular sequence. Under
$\Ext^1_R(I,R)\simeq R/I$, let $\bar r$ represent an extension
$0\to R\to M\to I\to0$.
Then $M$ is free at $\mathfrak m$ if and only if $\bar r$ is a unit in
$R/I$. If $\bar r$ is nonzero but vanishes only at a codimension-three
point, then $M$ has projective dimension one and depth two at that point.
\end{lemma}
\begin{proof}
Pushing out the Koszul resolution of $I$ by a lift $r$ of
$\bar r$ gives, up to signs,
\[
0\to R\xrightarrow{(-g,f,r)}R^3\to M\to0.
\]
The module $M$ is torsion-free of rank two, since it is an extension
of $I$ by $R$. Over the local ring $R$, the quotient $M$ is free
exactly when the inclusion splits, equivalently when its entries
generate the unit ideal.
\[
M\text{ free}\ \Longleftrightarrow\ (f,g,r)=R
\ \Longleftrightarrow\ \bar r\in(R/I)^\times.
\]
Otherwise the resolution is minimal, so $\operatorname{pd}_R M=1$.
The Auslander--Buchsbaum formula then gives
$\operatorname{depth}_R M=3-1=2$.
After localization, the same criterion identifies the nonfree locus
with $V(f,g,r)$, proving the assertion at the isolated zero.
\end{proof}

\begin{lemma}\label{orig:4.4}
Let $\sigma$ be a nonzero class in the extension space of Lemma~\ref{orig:4.2}, and let
$F_\sigma$ be the middle term of
\[
0\longrightarrow\OO_X(bA)\longrightarrow F_\sigma\longrightarrow
\II_Z(\zeta-bA)\longrightarrow0.\tag{3.4}\label{eq:orig:4.4}
\]
Then $F_\sigma$ is a rank-two reflexive sheaf. Its singular scheme,
defined by $\Fitt_2(F_\sigma)$, equals the zero scheme of
$\sigma\in H^0(C,\OO_C((2b+1)H))$.
It is locally free if and only if $\sigma$ is nowhere zero; since
$\deg\OO_C((2b+1)H)>0$, no nonzero extension in this family is locally
free. For a general class the singular scheme is reduced and has
length $(2b+1)d$.
\end{lemma}
\begin{proof}
The extension makes $F_\sigma$ torsion-free of rank two and locally
free off $Z$. For any $p\in Z$, put $R=\OO_{X,p}$, write
$\II_{Z,p}=(f,g)$, and choose local trivializations of the line
bundles. Under the identification of Lemma~\ref{orig:4.2}, the germ
$\sigma_p$ is represented by the local extension class
$\bar r\in R/(f,g)=\OO_{Z,p}$. Lemma~\ref{orig:4.3}, after
localization if necessary, gives
\[
(F_\sigma)_p\text{ free}
\quad\Longleftrightarrow\quad
\bar r\in\OO_{Z,p}^{\times}
\quad\Longleftrightarrow\quad
\sigma(p)\ne0.
\]
Indeed, $\bar r$ is a unit exactly when $\sigma_p$ generates the
line-bundle stalk, or equivalently when its image $\sigma(p)$ in
$\OO_Z((2b+1)A)\otimes\kappa(p)$ is nonzero.
Since $\sigma\ne0$ on the smooth integral curve $Z\simeq C$, its
zeros are finite. Lemma~\ref{orig:4.3} gives depth two at each zero;
elsewhere $F_\sigma$ is locally free. It satisfies
$S_2$ and is reflexive by the torsion-free $S_2$ criterion
\cite[\S1]{Hartshorne}.

For a lift $r\in R$ of $\bar r$, the relation matrix
$(-g,f,r)^{\mathsf t}$ also gives
\[
\Fitt_2((F_\sigma)_p)=(f,g,r),\qquad
\OO_{\Sing(F_\sigma),p}\simeq\OO_{Z,p}/(\sigma).
\]
The singularity scheme is $(\sigma)_0$, including multiplicities.

Every nonzero section has zero divisor of degree
$\deg\OO_C((2b+1)H)=(2b+1)d>0$, so no such extension is locally
free. The line bundle is globally generated because $H$ is
base-point-free. Bertini's theorem on $C$ gives a reduced zero
divisor for a general section \cite[III, Corollary 10.9]{HartshorneAG}.
\end{proof}

It remains to compute the Chern classes.

\begin{lemma}\label{orig:4.6}
The sheaf $F_\sigma$ has
\[
c_1(F_\sigma)=\zeta,\quad c_2(F_\sigma)=(b+1)\zeta A-b^2A^2,
\quad c_3(F_\sigma)=(2b+1)d.
\]
\end{lemma}
\begin{proof}
Put $g_0=bA$, $m_0=\zeta-bA$, $z=\zeta A$. The Koszul resolution gives
$\ch(\II_Z)=1-z+\tfrac12z(\zeta+A)$ through degree three. Thus
$\ch(F_\sigma)=\exp(g_0)+\exp(m_0)\ch(\II_Z)$ yields
\[
c_1=g_0+m_0=\zeta,\qquad c_2=g_0m_0+z=(b+1)\zeta A-b^2A^2.
\]
For the third class, compare with
$\mathcal E=\OO_X(g_0)\oplus\OO_X(m_0)$, for which
$c_1(\mathcal E)=c_1$, $c_2(\mathcal E)=g_0m_0$, and $c_3(\mathcal E)=0$.
Using $\ch_3=(c_1^3-3c_1c_2+3c_3)/6$, the difference is
\[
\frac{c_3(F_\sigma)-c_1z}{2}
=\ch_3(F_\sigma)-\ch_3(\mathcal E)
=-zm_0+\frac12z(\zeta+A).
 \]
Solving for the third Chern class gives
\[
c_3(F_\sigma)=z(\zeta+A+c_1-2m_0)
=z(\zeta+A+2g_0-c_1)=(2b+1)\zeta A^2=(2b+1)d.
\]
\end{proof}

The preceding lemmas identify the extension space, describe the
singularities and compute the Chern classes. We collect these properties
in the following theorem, which will serve as the starting point for the
stability and deformation arguments.

\begin{theorem}\label{orig:4.1}\label{thm:serre-ext}
Let $b\geq1$. For general $s\in H^0(S,\VV)$ and $C\in|H|$,
the complete intersection $Z=D_s\cap\pi^{-1}(C)$ is isomorphic to $C$, and
\[
\Ext^1_X(\II_Z(\zeta-bA),\OO_X(bA))
\simeq H^0(C,\OO_C((2b+1)H)).
\]
Every nonzero class $\sigma$ gives a rank-two reflexive extension
$0\to\OO_X(bA)\to F_\sigma\to\II_Z(\zeta-bA)\to0$ with
\[
(c_1,c_2,c_3)=\bigl(\zeta,(b+1)\zeta A-b^2A^2,(2b+1)d\bigr).
\]
Its Fitting singularity scheme is $(\sigma)_0$, of length $(2b+1)d$;
for general $\sigma$ it is reduced.
\end{theorem}
\begin{proof}
The identification $Z\simeq C$ is \eqref{eq:orig:4.2}.
Lemma~\ref{orig:4.2} identifies the extension space, and
Lemma~\ref{orig:4.4} proves reflexivity and describes the Fitting
singularity scheme. The Chern classes are computed in
Lemma~\ref{orig:4.6}.
\end{proof}

\subsection{Sections and extensions for curve quotients}
The enlarged Serre family and the boundary construction use
quotients on a smooth curve $C$ with kernels of degree zero;
the quotients may have torsion.

\begin{proposition}\label{orig:6.12}
Let $i:C\hookrightarrow S$ be a smooth member of $|H|$, and let
$\VV|_C\twoheadrightarrow J$ be a generically rank-one quotient,
possibly with torsion, whose kernel $K$ has degree zero.
Then $K$ is a line bundle. Put $Z=\PP_C(J)\subset X$ and let $b\ge1$.
\begin{enumerate}[label=\textup{(\roman*)}]
\item Every extension
\[
0\to\OO_X(bA)\to F\to\II_Z(\zeta-bA)\to0
\]
satisfies
\[
h^0(X,F)=b^2h+2+h^0(S,\VV(-(b+1)H)).
\tag{3.5}\label{eq:orig:6.8}
\]
\item Set $\mathcal L_K:=K^{-1}\otimes\OO_C((2b+1)H)$. Then
\[
\begin{aligned}
\Ext^1_X(\II_Z(\zeta-bA),\OO_X(bA))
&\simeq H^0(C,\mathcal L_K),\\
h^0(C,\mathcal L_K)&=R_b=(4b+1)h,\\
\Ext^2_X(\II_Z(\zeta-bA),\OO_X(bA))&=0.
\end{aligned}\tag{3.6}\label{eq:integral-extension}
\]
In fact, all Ext groups in degrees other than one vanish.
\item Every such extension has the numerical Chern classes of
Lemma~\ref{orig:4.6}:
\[
c_1(F)=\zeta,\qquad c_2(F)\equiv(b+1)\zeta A-b^2A^2,
\qquad \deg c_3(F)=(2b+1)d.
\]
\end{enumerate}
\end{proposition}
\begin{proof}
The kernel $K$ is a torsion-free rank-one sheaf on the smooth curve $C$,
hence a line bundle. Set $X_C=\pi^{-1}(C)$, with $p:X_C\to C$ and
$j:X_C\hookrightarrow X$. The kernel of
$\operatorname{Sym}(\VV|_C)\twoheadrightarrow\operatorname{Sym}J$
is generated by $K$ in degree one. The resulting nonzero map of line
bundles $p^*K\to\OO_{X_C}(\zeta)$ is injective on the integral surface $X_C$,
so $\II_{Z/X_C}(\zeta)\simeq p^*K$, even when $J$ has torsion.
The ideal sequence for $Z\subset X_C\subset X$ gives
\[
0\to\OO_X(\zeta-(b+1)A)\to\II_Z(\zeta-bA)
\to j_*p^*(K(-bH))\to0.
\]
Since $\deg K(-bH)=-2bh<0$ and $p_*\OO_{X_C}=\OO_C$,
the right-hand term has no global sections, so
\[
H^0(X,\II_Z(\zeta-bA))
\simeq H^0(S,\VV(-(b+1)H)).
\]
Now $H^1(X,\OO_X(bA))=0$ and $h^0(X,\OO_X(bA))=b^2h+2$
by Lemma~\ref{orig:2.8} and the cohomology vanishing on $S$, proving \textup{(i)}.

For \textup{(ii)}, put $\mathcal Q_Z=\II_Z(\zeta-bA)$.
The line-bundle term of the ideal sequence contributes no Ext, since
\[
R\Hom_X(\OO_X(\zeta-(b+1)A),\OO_X(bA))
\simeq R\Gamma(X,\OO_X(-\zeta+(2b+1)A))=0
\]
by Lemma~\ref{orig:2.8}. Cartier-divisor duality
\cite[III, \S\S6--7]{HartshorneRD} gives
$j^!\OO_X(bA)\simeq p^*\OO_C((b+1)H)[-1]$, and $Rp_*\OO_{X_C}=\OO_C$.
\[
\begin{aligned}
R\Hom_X(\mathcal Q_Z,\OO_X(bA))
&\simeq R\Hom_{X_C}\bigl(p^*(K(-bH)),p^*\OO_C((b+1)H)\bigr)[-1]\\
&\simeq R\Gamma(C,\mathcal L_K)[-1].
\end{aligned}
\]
This identifies $\Ext^1$ as claimed. Since $\omega_C=\OO_C(H)$,
Serre duality gives
\[
\Ext^2_X(\mathcal Q_Z,\OO_X(bA))
\simeq H^1(C,\mathcal L_K)
\simeq H^0(C,K(-2bH))^\vee=0.
\]
Finally, $g(C)=h+1$ and $\deg\mathcal L_K=(4b+2)h$,
so Riemann--Roch gives $\dim\Ext^1_X(\mathcal Q_Z,\OO_X(bA))=(4b+1)h$.
The derived identification also gives vanishing in every other degree.

For \textup{(iii)}, we use the same ideal sequence to compare numerical
Chern characters. Since $\deg K=0$, $p^*K$ is numerically trivial on
$X_C$. Grothendieck--Riemann--Roch for $j:X_C\hookrightarrow X$ gives
\[
\ch\bigl(j_*p^*(K(-bH))\bigr)
\equiv\ch\bigl(j_*\OO_{X_C}(-bA)\bigr)
=e^{-bA}-e^{-(b+1)A},
\]
where the last equality follows from the Cartier divisor sequence for
$X_C$. Additivity in the ideal and extension sequences therefore yields
\[
\ch(F)\equiv e^{bA}+e^{\zeta-(b+1)A}+e^{-bA}-e^{-(b+1)A}.
\]
This is the character computed from the Koszul resolution in
Lemma~\ref{orig:4.6}, so its numerical Chern classes are the stated ones.
The equality for $c_1$ also holds as a divisor class, directly from
the determinant of the extension.
\end{proof}

\begin{remark}\label{rem:integral-quotients}
No splitting assumption on $\VV$ or $\pi_*F$ is needed.
The quotients $\VV|_C\twoheadrightarrow(\det\VV)|_C$ in the
construction of Theorem~\ref{thm:serre-ext} have kernel $\OO_C$
and satisfy the hypotheses above. Proposition~\ref{orig:6.13}
constructs larger families of line-bundle quotients with
kernels of degree zero.
\end{remark}

\section{Stability intervals and numerical walls}\label{sec:stability}
\subsection{Stability for arbitrary Picard rank}
We first prove stability for the polarizations satisfying
$2bAL^2<\zeta L^2<2(b+1)AL^2$.
Under \eqref{eq:N}, these polarizations form $\CC_b$, and we then show
that no numerically effective rank-one wall meets $\CC_b$.

For the subsheaf analysis we isolate a numerical condition on the
polarization:
\[\label{eq:slope-window}
\pi^*c_1(\mathcal V)\cdot L^2<2bA\cdot L^2<\zeta\cdot L^2
<2(b+1)A\cdot L^2. \tag{$\dagger$}
\]

\begin{lemma}\label{lem:subsheaf-test}
Let $b\ge1$, let $L\in\mathcal A_+$ satisfy~\eqref{eq:slope-window}, and let
\begin{equation*}
0\to\mathcal O_X(bA)\to F\to\mathcal I_Z(\zeta-bA)\to0
\end{equation*}
be a nonzero torsion-free extension with $Z\ne\emptyset$ of
codimension two. Assume one of the following:
\begin{enumerate}[label=\textup{(\alph*)}]
\item $\operatorname{Pic}(S)=\Z H$;
\item $F$ is reflexive, and $Z\subset\pi^{-1}(C)$ maps onto a
smooth curve $C\in|H|$.
\end{enumerate}
Then $F$ is $L$-slope stable.
\end{lemma}

\begin{proof}
Put $\mathcal Q_Z=\mathcal I_Z(\zeta-bA)$, so that
$\mu_L(F)=\tfrac12\zeta\cdot L^2$. Let $M_1\subset F$ be a saturated
subsheaf of rank one. Then $M_1\simeq\mathcal I_T\otimes L_1$ with $L_1$ a
line bundle and $\operatorname{codim}T\ge2$; if $F$ is reflexive, so is
$M_1$, and $T=\emptyset$. Write $c_1(L_1)=a\zeta+\pi^*\beta_1$ with
$\beta_1\in\operatorname{Pic}(S)$.

For a general point $x\in S$, the fiber $f=\pi^{-1}(x)$ is disjoint
from $Z$ and $T$. By generic flatness, $F/M_1$ is flat over a dense open subset of $S$, so
the restriction $M_1|_f\to F|_f$ is injective for general $f$. Since
$\operatorname{Ext}^1_{\PP^1}(\mathcal O(1),\mathcal O)=0$, we have
$F|_f\simeq\mathcal O_{\PP^1}\oplus\mathcal O_{\PP^1}(1)$,
and $M_1|_f\simeq\mathcal O_{\PP^1}(a)$ is a subsheaf of it. Hence
$a\le1$.

Let $\varphi\colon M_1\to \mathcal Q_Z$ be the induced map. If $\varphi=0$, then
$M_1\subset \OO_X(bA)$ and
$\mu_L(M_1)\le bA\cdot L^2<\mu_L(F)$ by~\eqref{eq:slope-window}.

Assume $\varphi\ne0$. If $M_1\cap \OO_X(bA)\neq0$, it has rank one, and the
torsion sheaf $M_1/(M_1\cap \OO_X(bA))$ injects into the torsion-free $\mathcal Q_Z$, forcing
$\varphi=0$. Hence $\varphi$ is injective, and its double dual
$\varphi^{\vee\vee}\colon L_1\to\mathcal O_X(\zeta-bA)$ is a nonzero map
of line bundles; consequently, $\mathcal O_X((1-a)\zeta)\otimes\pi^*\mathcal O_S(-bH-\beta_1)$ admits a nonzero global section.

If $a\le0$, pushing forward gives
$H^0(S,\operatorname{Sym}^{1-a}\mathcal V\otimes\mathcal O_S(-bH-\beta_1))\ne0$.
Since $\mathcal V$ is globally generated,
$\mathcal V\simeq\mathcal V^\vee\otimes\det\mathcal V
\hookrightarrow(\det\mathcal V)^{\oplus n}$ for some $n$, hence
$\operatorname{Sym}^{1-a}\mathcal V\hookrightarrow
\bigl((\det\mathcal V)^{\otimes(1-a)}\bigr)^{\oplus r_a}$. Here $r_a$ is a positive integer. A nonzero component of the
section shows that $D_1:=(1-a)c_1(\mathcal V)-bH-\beta_1$ is effective.
\[
\mu_L(M_1)=a\,(\zeta-\pi^*c_1(\mathcal V))\cdot L^2
+(\pi^*c_1(\mathcal V)-bA)\cdot L^2-\pi^*D_1\cdot L^2
\le(\pi^*c_1(\mathcal V)-bA)\cdot L^2,
\]
because $a\le0$, $(\zeta-\pi^*c_1(\mathcal V))\cdot L^2>0$
by~\eqref{eq:slope-window}, and $\pi^*D_1\cdot L^2\ge0$. By
\eqref{eq:slope-window} the right-hand side is less than
$bA\cdot L^2<\mu_L(F)$.

If $a=1$, then $\varphi^{\vee\vee}=\pi^*u$ for a nonzero
$u\in H^0(S,\mathcal O_S(-bH-\beta_1))$. Let $D_2$ be its zero divisor, so that
$c_1(L_1)=\zeta-bA-\pi^*D_2$.

Under (b), $T=\emptyset$, so $\pi^*u$ vanishes on $Z$. Since $Z$ maps
onto $C$, $u$ vanishes on $C$ and $D_2-H$ is effective. Hence
$\mu_L(M_1)\le(\zeta-(b+1)A)\cdot L^2<\mu_L(F)$
by~\eqref{eq:slope-window}.

Under (a), $D_2\equiv mH$ with $m\ge0$. If $m\ge1$, the same bound
applies. If $m=0$, the effective divisor $D_2$ is numerically
trivial, hence zero. Thus $u$ is nowhere zero and
$\varphi^{\vee\vee}$ is an isomorphism. Outside $T$, the inclusion
$M_1\hookrightarrow\mathcal Q_Z$ is therefore an isomorphism, so
$\mathcal G'=\mathcal Q_Z/M_1$ is supported on $T$, in codimension at
least two. Hence $\shExt^i_X(\mathcal G',\OO_X(bA))=0$ for $i\le1$
\cite[Proposition~1.1.6(i)]{HL}, and the local-to-global spectral sequence
gives $\operatorname{Ext}^i_X(\mathcal G',\OO_X(bA))=0$ for $i\le1$. Applying
$\operatorname{Hom}_X(-,\OO_X(bA))$ to $0\to M_1\to \mathcal Q_Z\to\mathcal G'\to0$ shows
that
\[
\operatorname{Ext}^1_X(\mathcal Q_Z,\OO_X(bA))\to\operatorname{Ext}^1_X(M_1,\OO_X(bA))
\]
is injective. The image of $[F]$ is the class of the pullback of the
extension along $\varphi$, which is split by $M_1\subset F$, so $[F]=0$, a contradiction.

In all cases $\mu_L(M_1)<\mu_L(F)$.
\end{proof}

\begin{theorem}\label{orig:5.2}\label{thm:arbitrary-stability}
Retain the standing assumptions, without assuming \eqref{eq:N}, and let
$b\ge1$ satisfy $2bH^2>c_1(\VV)\cdot H$.
The polarizations $L\in\mathcal A_+$ satisfying
\[
2b\,A\cdot L^2<\zeta\cdot L^2<2(b+1)\,A\cdot L^2
\]
form a nonempty open subset of $\mathcal A_+$.
Every nonzero extension of Theorem~\ref{thm:serre-ext} is slope
stable for every such $L$. For rational $L$, the Brill--Noether loci $W^k_{X,L}(w)$
are nonempty for $1\le k\le b^2h+2$, where $w$ is the common
numerical Chern character of these extensions.
\end{theorem}

\begin{proof}
Since stability is unchanged by positive scaling, write
$L_t=\zeta+tA$, with $t>0$. To determine the interval of $t$
satisfying the stated inequalities, we use global generation
and the Hodge index theorem:
\[
c_1(\VV)H\ge0,\qquad\gamma\ge0,\qquad
c_1(\VV)^2\le\frac{(c_1(\VV)H)^2}{d}.
\]
For $m=b,b+1$, Lemma~\ref{orig:2.7} gives
\[
\begin{aligned}
f_m(t)&=(\zeta-2mA)\cdot L_t^2\\
&=dt^2+2(c_1(\VV)H-2md)t+c_1(\VV)^2-\gamma-2mc_1(\VV)H.
\end{aligned}
\]
The hypothesis $2bd>c_1(\VV)H$ makes the linear coefficient
negative, while
\[
f_m(0)\le c_1(\VV)H\left(\frac{c_1(\VV)H}{d}-2m\right)-\gamma\le0.
\]
Hence $f_m$ has a unique positive root $t_m$: if $f_m(0)=0$,
the other root is positive because the linear coefficient is negative.
Moreover,
\[
f_{b+1}(t_b)=-2A\cdot L_{t_b}^2<0,
\]
so $t_b<t_{b+1}$. The required inequalities are
$f_b(t)>0>f_{b+1}(t)$, and hold exactly on
the nonempty interval $(t_b,t_{b+1})$.

To control subsheaves of nonpositive degree on the fibers, we need
one further intersection bound. For every $t>0$, the same Hodge
inequality gives directly
\[
\begin{aligned}
(2bA-\pi^*c_1(\VV))\cdot L_t^2
&=2b\,c_1(\VV)H-c_1(\VV)^2+2t(2bd-c_1(\VV)H)\\
&\ge(2bd-c_1(\VV)H)
\left(2t+\frac{c_1(\VV)H}{d}\right)>0.
\end{aligned}
\]
Thus every $L=L_t$ with $t_b<t<t_{b+1}$ satisfies
\[
(\pi^*c_1(\VV))\cdot L^2<2bA\cdot L^2<\zeta\cdot L^2.
\]

Together with $f_{b+1}(t)<0$, that is, $\zeta\cdot L^2<2(b+1)A\cdot L^2$,
every such $L$ satisfies~\eqref{eq:slope-window}. Every nonzero
extension $F$ of Theorem~\ref{thm:serre-ext} is reflexive by Lemma~\ref{orig:4.4}, and its curve
$Z\simeq C$ lies in $\pi^{-1}(C)$ and maps isomorphically onto the
smooth curve $C\in|H|$. Lemma~\ref{lem:subsheaf-test}(b) shows
that $F$ is slope stable.

Finally, the distinguished line subsheaf gives
$h^0(F)\ge h^0(S,\OO_S(bH))=b^2h+2$.
For rational $L$, slope stability implies Gieseker stability,
so $[F]\in W^k_{X,L}(w)$ for $1\le k\le b^2h+2$.
\end{proof}

\subsection{Absence of numerical walls for Picard rank one}
The line subsheaf $\OO_X(bA)$ determines a slope-equality wall.
To describe the polarizations of Theorem~\ref{thm:arbitrary-stability}
in terms of numerical walls,
we now assume \eqref{eq:N} and fix the following notation.

\begin{definition}\label{orig:3.1}\label{def:stability-region}
Put $c_1(\VV)\equiv eH$, and fix an integer $b\geq1$.
For the extension construction, set
\[
c_1=\zeta,\qquad c_2=(b+1)\zeta A-b^2A^2,\qquad
\xi_b=-\zeta+2bA.\tag{4.1}\label{eq:orig:3.1}
\]
The associated \emph{distinguished wall} is $W^{\xi_b}$.
For $b>e/2$, define
\[
\CC_b=\{L\in\mathcal A_+:\xi_bL^2<0<\xi_{b+1}L^2\}.
\]
\end{definition}

For the fixed Chern classes \eqref{eq:orig:3.1}, one has
\[
\Delta_{\xi_b}=\zeta A,\qquad
\Delta_{\xi_{b+1}}=(2b+1)A^2.\tag{4.2}\label{eq:orig:3.2}
\]
Both classes are effective: $\zeta A$ is represented by the
complete-intersection curve of Theorem~\ref{thm:serre-ext},
and $A^2$ is a sum of fibers.
With $t_b<t_{b+1}$ denoting the positive roots defined in the proof
of Theorem~\ref{thm:arbitrary-stability}, the set
\[
\CC_b=\{\alpha(\zeta+tA):\alpha>0,\ t_b<t<t_{b+1}\}
\tag{4.3}\label{eq:orig:5.1}
\]
is nonempty. In particular every $L\in\mathcal C_b$ satisfies~\eqref{eq:slope-window},
since $\pi^*c_1(\mathcal V)\cdot L^2=eA\cdot L^2<2bA\cdot L^2$ when
$2b>e$. By Lemma~\ref{lem:subsheaf-test}(a), every nonzero torsion-free
extension $0\to\OO_X(bA)\to F\to\II_Z(\zeta-bA)\to0$ with
$Z\ne\varnothing$ and Chern classes \eqref{eq:orig:3.1} is slope-stable for every $L\in\CC_b$. The next proposition shows that no numerically effective rank-one wall meets $\CC_b$, which is needed for Corollary~\ref{cor:chamber}.

\begin{proposition}\label{orig:5.3}\label{thm:rank-one-stability}
Assume \eqref{eq:N}, $b\geq1$, and $2b>e$. For the Chern classes
\eqref{eq:orig:3.1}, the set $\CC_b$ contains no numerically effective
rank-one wall and lies in a chamber adjacent to $W^{\xi_b}$.
\end{proposition}
\begin{proof}
To rule out numerically effective walls in $\CC_b$, we intersect the effective residual class of a putative wall with the
nef divisor $A$.

Up to sign every wall class has the form
\[
\eta=-(2a+1)\zeta+2pA,\qquad a\geq0,\quad p\in\Z.\tag{4.4}\label{eq:orig:5.2}
\]
Using $\zeta^2\equiv e\zeta A-(\gamma/d)A^2$ gives
\[
\begin{gathered}
\Delta_\eta=\alpha_\eta\zeta A+\beta_\eta A^2,\\
\alpha_\eta=b+1+a(a+1)e-(2a+1)p,\qquad
\beta_\eta=p^2-b^2-a(a+1)\gamma/d .
\end{gathered}\tag{4.5}\label{eq:orig:5.3}
\]
Effectivity and nefness of $A$ imply
\[
A\Delta_\eta=d\alpha_\eta\geq0.\tag{4.6}\label{eq:orig:5.4}
\]
Put $n=2a+1$. If $W^\eta$ met $\CC_b$, the identities
$\eta=n\xi_b+2(p-nb)A=n\xi_{b+1}+2(p-n(b+1))A$
would give
\[
nb<p<n(b+1).\tag{4.7}\label{eq:orig:5.5}
\]
For $a=0$ this contains no integer. For $a\geq1$, $p\geq nb+1$
gives $\alpha_\eta\leq a(a+1)(e-4b)-2a<0$, a contradiction.
Thus no numerically effective rank-one wall meets $\CC_b$.
The interval description \eqref{eq:orig:5.1} implies that $\CC_b$ is connected, and its boundary wall $W^{\xi_b}$ is numerically effective
by \eqref{eq:orig:3.2}. This gives the chamber assertion.

\end{proof}

\begin{corollary}\label{cor:chamber}
Assume \eqref{eq:N}, $b\ge1$, and $2b>e$. Fix a numerical Chern
character $w$ of rank two with first two Chern classes
\eqref{eq:orig:3.1}. For every $L\in\CC_b$, a torsion-free sheaf
of character $w$ is slope-semistable if and only if it is
slope-stable, and this condition is independent of $L\in\CC_b$.
For rational $L\in\CC_b$, Gieseker stability coincides
with slope stability, and the moduli spaces $M_{X,L}(w)$ are
canonically isomorphic as $L$ varies in $\CC_b$.
\end{corollary}
\begin{proof}
Let $F$ have character $w$, and let $M_1\subset F$ be a saturated
rank-one subsheaf. Set $\eta=2c_1(M_1)-\zeta$. The class $\eta$ is
nonzero, since the coefficient of $\zeta$ is odd in
$\Pic(X)=\Z\zeta\oplus\Z A$. The computation following
Definition~\ref{def:walls} shows that $\Delta_\eta$ is effective.
The product $\eta\cdot L^2$ never vanishes on $\CC_b$, by
Proposition~\ref{thm:rank-one-stability}. Since $\CC_b$ is connected by \eqref{eq:orig:5.1}, this continuous
function has constant nonzero sign there. Applying this to every
$M_1$ proves the assertions about slope stability.
The implications slope-stable $\Rightarrow$ Gieseker-stable
$\Rightarrow$ slope-semistable give the claimed equivalence.
The stable families of character $w$ are the same for all rational
$L\in\CC_b$; their coarse moduli spaces are canonically
isomorphic.
\end{proof}

\section{Moduli components and tangent corrections}\label{sec:components}
Fix $b\ge1$ and a rational polarization $L$ in the stability range
of Theorem~\ref{thm:arbitrary-stability}. Let $w=\ch(F_\sigma)$ be
the numerical Chern character of the Serre sheaves
\eqref{eq:orig:4.4}, and write
\[
\FF_b=f_b(\mathcal T_b)\subset M_{X,L}(w)
\]
for the image of the extension family in Lemma~\ref{orig:6.1}. Put
\[
D_b=\dim\FF_b,\qquad
\varepsilon_F=\dim\Ext^1_X(F,F)-D_b.
\]
Put $M=M_{X,L}(w)$, and let $\mathcal M_b$ be the closure
$\overline{\FF_b}$ in $M$, endowed with its reduced structure.
We compute $\varepsilon_F$ and use its vanishing to characterize
generically smooth moduli components.

\subsection{The extension family and component criterion}\label{subsec:extension-family}
Before treating the individual parameter spaces, we record the extension
construction that will be used for both line-bundle and torsion quotients.

\begin{lemma}\label{lem:serre-family}
Let $b\ge1$ satisfy $2bd>c_1(\VV)\cdot H$, and let $L$ be a rational
polarization. Let $B$ be a smooth irreducible variety and
$\ZZ\subset X\times B$ a $B$-flat family of codimension-two local
complete intersection curves. Suppose that each geometric fiber is
$Z_y=\PP_{C_y}(J_y)$ for a quotient as in
Proposition~\ref{orig:6.12}, that $y\mapsto Z_y$ is injective on
geometric points, and that every nonzero extension of
$\II_{Z_y}(\zeta-bA)$ by $\OO_X(bA)$ is $L$-slope stable.
\begin{enumerate}[label=\textup{(\roman*)}]
\item The relative extension sheaf $\EE_B$ of
$\mathcal Q_B=\II_{\ZZ}\otimes\pr_X^*\OO_X(\zeta-bA)$ by
$\mathcal{L}_B=\pr_X^*\OO_X(bA)$ is locally free of rank $R_b$ and
commutes with arbitrary base change.
\item On $T=\PP_B(\EE_B)$, parametrizing lines of extension classes,
there is a flat tautological extension
\[
0\to\mathcal L_T\otimes p_T^*\OO_T(1)
\to\widetilde{\FF}\to\mathcal Q_T\to0,
\]
where $p_T:X\times T\to T$, and $\mathcal L_T,\mathcal Q_T$ are the
pullbacks of $\mathcal{L}_B,\mathcal Q_B$. Its fibers have numerical
character $w$ and define a
classifying morphism $f:T\to M_{X,L}(w)$.
\item The morphism $f$ is injective on geometric points. Its image
is irreducible of dimension $\dim B+R_b-1$.
\end{enumerate}
\end{lemma}
\begin{proof}
The flatness of $\ZZ$ makes $\mathcal Q_B$ flat over $B$, and it is
perfect because $X\times B$ is smooth. For $p_B:X\times B\to B$,
proper perfect pushforward and derived base change
\cite[Theorem~2.5.4]{TT}, \cite[Exp.~IV, Proposition~3.1.0]{SGA6}
apply to $Rp_{B*}R\shHom(\mathcal Q_B,\mathcal{L}_B)$.
By Proposition~\ref{orig:6.12}\textup{(ii)}, every derived fiber is
concentrated in degree one, of dimension $R_b$. Locally, removing
acyclic summands from a finite free complex leaves only its degree-one
term \cite[Lemma~2.1.3]{Lieblich}. Thus
\[
Rp_{B*}R\shHom(\mathcal Q_B,\mathcal{L}_B)\simeq\EE_B[-1],
\qquad \rk\EE_B=R_b,
\]
and the same description after arbitrary base change proves \textup{(i)}.

To assemble the extensions globally, let $\Pi:T\to B$. The
tautological inclusion $\OO_T(-1)\hookrightarrow\Pi^*\EE_B$ is a
section of $\Pi^*\EE_B\otimes\OO_T(1)$. Derived base change and the
local-to-global spectral sequence on $T$, using concentration in degree
one, identify
\[
\Ext^1_{X\times T}(\mathcal Q_T,\mathcal L_T\otimes p_T^*\OO_T(1))
\simeq H^0(T,\Pi^*\EE_B\otimes\OO_T(1)).
\]
The tautological section therefore defines the asserted extension,
whose fiber classes are nonzero. Both end terms are $T$-flat, hence
so is the middle term. Proposition~\ref{orig:6.12}\textup{(iii)}
gives its numerical character $w$, and the assumed stability gives $f$.

It remains to recover the parameters from a fiber $F$. The inclusion
$\VV\hookrightarrow(\det\VV)^{\oplus h^0(S,\VV)}$ and
$2bd>c_1(\VV)\cdot H$ give $H^0(S,\VV(-2bH))=0$: a line bundle
of negative intersection with the ample class $H$ has no section.
Twisting the extension by $-bA$ now gives
\[
H^0(X,\II_Z(\zeta-2bA))\subset H^0(S,\VV(-2bH))=0,
\qquad h^0(F(-bA))=1.
\]
Thus $F$ determines $\OO_X(bA)\hookrightarrow F$ up to scalar.
An isomorphism of its quotients
$\II_Z(\zeta-bA)\simeq\II_{Z'}(\zeta-bA)$ extends to their double
duals as a scalar automorphism of $\OO_X(\zeta-bA)$; naturality of
the double-dual inclusions gives $\II_Z=\II_{Z'}$.
The hypothesis on the curve family then recovers $y\in B$.
Both end terms have only scalar automorphisms, so the extension class
is determined up to scalar as well. This proves injectivity on
geometric points, hence quasi-finiteness. Since $T$ is irreducible,
its image has dimension $\dim T=\dim B+R_b-1$.
\end{proof}

For the remainder of this subsection, we assume \eqref{eq:N}, $b\ge1$, and $2b>e$,
and let $L\in\CC_b$ be rational.
Proposition~\ref{prop:5.7} gives the large-$b$ analogue without \eqref{eq:N}.

Since $b\ge1$, Kodaira vanishing on the K3 base gives
\[
H^i(\OO_S(2bH))=H^i(\OO_S((2b+1)H))=0\qquad(i>0).
\]
The inclusion $\VV\hookrightarrow\OO_S(eH)^{\oplus h^0(S,\VV)}$ and
$e-2b<0$ also give $H^0(\VV(-2bH))=0$.
We use these vanishings throughout this subsection.

Let $U_0\subset\PP H^0(S,\VV)\times|H|$ be the nonempty open
subset of pairs $([s],C)$ such that $C$ is smooth and $s|_C$ is
nowhere zero. Define
\[
\mathscr C_{U_0}
=\{(x,([s],C))\in S\times U_0:x\in C\}.
\]
On $(\pi\times\mathrm{id}_{U_0})^{-1}(\mathscr C_{U_0})$,
the universal section $s$, viewed as a section of the tautological
line bundle with its parameter twist, defines a closed subscheme
\[
\ZZ_{U_0}\subset X\times U_0,\qquad
(\ZZ_{U_0})_{([s],C)}=(s=0)_X\cap\pi^{-1}(C).
\]
Since $s|_C$ is nowhere zero, the quotient-line construction in
\eqref{eq:orig:4.2}, applied in families, identifies $\ZZ_{U_0}$
with $\mathscr C_{U_0}$ over $U_0$.
The morphism $\ZZ_{U_0}\to U_0$ is smooth and projective, hence flat.
With respect to a fixed ample line bundle on $X$, its fibers have
a common Hilbert polynomial $P$, since $U_0$ is irreducible.
This family defines a morphism
\[
\eta_{\mathrm{Hilb}}:U_0\longrightarrow\Hilb^P(X),\qquad
([s],C)\longmapsto[(s=0)_X\cap\pi^{-1}(C)].
\]

For each $[Z]\in\eta_{\mathrm{Hilb}}(U_0)$, the fiber $\eta_{\mathrm{Hilb}}^{-1}([Z])$ identifies
with the nonempty open subset of $\PP H^0(X,\II_Z(\zeta))$
whose corresponding sections of $\VV$ are nowhere zero on
$C=\pi(Z)$. The Koszul resolution of $Z$, twisted by $\zeta$,
and $H^0(S,\OO_S(-H))=H^1(S,\OO_S(-H))=0$ give
\[
h^0(X,\II_Z(\zeta))=1+h^0(S,\VV(-H)).
\]
The fiber dimension is
\[
\dim\eta_{\mathrm{Hilb}}^{-1}([Z])
=h^0(X,\II_Z(\zeta))-1=h^0(S,\VV(-H)).
\]

Choose a smooth dense open subset $B$ of
$(\overline{\eta_{\mathrm{Hilb}}(U_0)})_{\mathrm{red}}$ contained in $\eta_{\mathrm{Hilb}}(U_0)$.
The fiber dimension formula gives
\[
\dim B=h+h^0(S,\VV)-h^0(S,\VV(-H)).
\]
The universal Hilbert family restricts to a smooth family
$\ZZ\subset X\times B$, whose inclusion is a regular immersion
of codimension two. We use this family in the extension construction
below; passage to this dense locus does not change the reduced
closure of the resulting Serre family.

\begin{lemma}\label{orig:6.1}
Assume \eqref{eq:N}. Fix $b\ge1$ with $2b>e$ and a rational $L\in\CC_b$.
There is an irreducible parameter variety $\mathcal T_b$, a $\mathcal T_b$-flat family
of stable reflexive sheaves on $X$, and a morphism
\[
f_b:\mathcal T_b\longrightarrow M_{X,L}(w)
\]
injective on geometric points. Its image $\FF_b=f_b(\mathcal T_b)$ has dimension
\[
D_b=(4b+2)h+h^0(S,\VV)-h^0(S,\VV(-H))-1.
\]
Every member belongs to $W^k_{X,L}$ for $k\leq h^0(S,\OO_S(bH))$.
\end{lemma}
\begin{proof}
The base $B$ constructed above is a smooth subvariety of the Hilbert
scheme, so $y\mapsto Z_y$ is injective on geometric points.
Theorem~\ref{thm:arbitrary-stability} gives stability for every nonzero
extension. Lemma~\ref{lem:serre-family} therefore supplies
$\mathcal T_b$, the flat family and the injective classifying morphism
$f_b$, with
\[
\dim\FF_b=\dim\mathcal T_b=\dim B+R_b-1=D_b.
\]
Reflexivity follows from Lemma~\ref{orig:4.4}.
Finally, the inclusion $\OO_X(bA)\hookrightarrow F$ yields the
stated BN range.
\end{proof}

A two-term locally free presentation of the extension allows us
to compute $\Ext^1_X(F,F)$.
\begin{lemma}\label{orig:6.3}
For $b\ge1$, every extension
$F_\sigma$ in \eqref{eq:orig:4.4} admits a resolution
\[
0\longrightarrow P\xrightarrow{\varphi}Q\longrightarrow F_\sigma
\longrightarrow0,\tag{5.1}\label{eq:orig:6.1}
\]
where
\[
P=\OO_X(-(b+1)A),\qquad
Q=\OO_X(bA)\oplus\OO_X(-bA)\oplus\OO_X(\zeta-(b+1)A).
\]
If $u\in H^0(X,\OO_X(A))$ and $s\in H^0(X,\OO_X(\zeta))$
define $\pi^{-1}(C)$ and $D_s$, respectively, then, with the displayed
order of the summands of $Q$ and up to Koszul signs,
\[
\varphi=(\tau,u,s),\qquad
\tau\in H^0(X,\OO_X((2b+1)A)),
\]
and $\tau|_Z$ represents $\sigma$.
\end{lemma}
\begin{proof}
Twist \eqref{eq:orig:4.3} by $\OO_X(\zeta-bA)$ and put
$Q'=\OO_X(-bA)\oplus\OO_X(\zeta-(b+1)A)$. This gives
\[
0\longrightarrow P\longrightarrow Q'\longrightarrow
\II_Z(\zeta-bA)\longrightarrow0.
\]
Kodaira vanishing and Lemma~\ref{orig:2.8} give
\[
\Ext^1_X(Q',\OO_X(bA))
=H^1(X,\OO_X(2bA))\oplus H^1(X,\OO_X(-\zeta+(2b+1)A))=0.
\]
Applying $\Hom_X(-,\OO_X(bA))$ therefore gives the exact sequence
\[
\Hom_X(Q',\OO_X(bA))\to\Hom_X(P,\OO_X(bA))
\to\Ext^1_X(\II_Z(\zeta-bA),\OO_X(bA))\to0.
\]
Every extension is the pushout by
$\tau\in\Hom(P,\OO_X(bA))=H^0(\OO_X((2b+1)A))$, giving
\eqref{eq:orig:6.1}. The remaining components have degrees
$A$ and $\zeta$, respectively, so $\varphi=(\tau,u,s)$.

Under the edge isomorphism of Lemma~\ref{orig:4.2}, the extension
class is represented by $\tau|_Z$. Surjectivity of restriction
also follows on $S$,
since the first term of the exact sequence
\[
0\to\OO_S(2bH)\xrightarrow{u}\OO_S((2b+1)H)
\to\OO_C((2b+1)H)\to0
\]
has vanishing $H^1$ by Kodaira vanishing. Changes of the Koszul signs
affect neither the cokernel's isomorphism class nor its Fitting ideal.
\end{proof}

Let $\FF_b$ be the family of Lemma~\ref{orig:6.1}, of dimension
\[
D_b=(4b+2)h+h^0(S,\VV)-h^0(S,\VV(-H))-1.
\tag{5.2}\label{eq:orig:6.2}
\]
At a stable member, $T_{[F]}M_{X,L}(w)\simeq\Ext^1_X(F,F)$
\cite{HL,VermeireModuli,Kleppe}.
The presentation \eqref{eq:orig:6.1} identifies the excess over
$D_b$ with the kernel of the map $\mu_F$ below. The criterion
concerns smoothness of the ambient moduli space along the component.

\begin{theorem}\label{orig:6.5}\label{thm:regular-component}
Assume \eqref{eq:N}, $b\ge1$ and $2b>e$, and let $L\in\CC_b$ be rational. For a stable sheaf $F=F_\sigma$ constructed in \eqref{eq:orig:4.4}, use the presentation $(\tau,u,s)$ of \eqref{eq:orig:6.1} and set
\[
\begin{aligned}
\mu_F:\ H^1(X,F(-bA))\oplus H^1(X,F(bA))
&\longrightarrow H^1(X,F((b+1)A)),\\
(\alpha,\beta)&\longmapsto\tau\alpha+u\beta.
\end{aligned}\tag{5.3}\label{eq:orig:6.3}
\]
Then
\[
\dim\Ext^1_X(F,F)=D_b+\dim\ker\mu_F.
\]
The following are equivalent:
\begin{enumerate}[label=\textup{(\roman*)}]
\item $\mu_F$ is injective for a general $[F]\in\FF_b$;
\item $\mathcal M_b$ is an irreducible component of $M_{\mathrm{red}}$, and
$M$ is smooth at a general point of $\mathcal M_b$.
\end{enumerate}
In this case $\dim \mathcal M_b=D_b$, and a general member has reduced
singularity scheme of length $(2b+1)d$.
A sufficient condition for \textup{(i)} is
$H^1(F(-bA))=H^1(F(bA))=0$ at a general member.
\end{theorem}
\begin{proof}
The vanishings at the start of this subsection give the
presentation \eqref{eq:orig:6.1}. To compute $\Ext^1_X(F,F)$,
apply $\Hom_X(-,F)$ to this two-term locally free resolution. The third summand
of $Q^\vee\otimes F$ has vanishing first cohomology:
twisting \eqref{eq:orig:4.4} gives
\[
H^1(X,F(-\zeta+(b+1)A))\simeq H^1(X,\II_Z(A))=0.
\]
Indeed the Koszul sequence twisted by $A$ is
\[
0\to\OO_X(-\zeta)\to\OO_X(A-\zeta)\oplus\OO_X
\to\II_Z(A)\to0.
\]
The terms of relative degree $-1$ are acyclic, so
$H^1(\II_Z(A))=H^1(\OO_X)=H^1(\OO_S)=0$ and
$h^0(\II_Z(A))=1$.
Thus $\Ext^1_X(Q,F)=H^1(F(-bA))\oplus H^1(F(bA))$,
and the map to $\Ext^1_X(P,F)=H^1(F((b+1)A))$ is $\mu_F$.

Stability gives $\Hom_X(F,F)=\C\cdot\mathrm{id}_F$ \cite{HL}.
The long exact sequence gives
\[
0\to\coker\bigl[\Hom_X(Q,F)\to\Hom_X(P,F)\bigr]
\to\Ext^1_X(F,F)\to\ker\mu_F\to0.
\]
Taking dimensions gives
\[
\begin{aligned}
\dim\Ext^1_X(F,F)
&=h^0(X,P^\vee\otimes F)-h^0(X,Q^\vee\otimes F)+1\\
&\hspace{2em}+\dim\ker\mu_F.
\end{aligned}
\tag{5.4}\label{eq:orig:6.4}
\]

Twisting \eqref{eq:orig:4.4}
by the summands of $P^\vee$ and $Q^\vee$ yields
\[
\begin{aligned}
h^0(X,F((b+1)A))&=\nu_{2b+1}+\nu_1+h^0(S,\VV)-1,\\
h^0(X,F(-bA))&=1,\\
h^0(X,F(bA))&=\nu_{2b}+1+h^0(S,\VV(-H)),\\
h^0(X,F(-\zeta+(b+1)A))&=1.
\end{aligned}
\]
Indeed, Kodaira vanishing gives
$H^1(X,\OO_X((2b+1)A))=H^1(X,\OO_X(2bA))=0$, so the
extension remains exact on global sections after twisting by
$(b+1)A$ and $bA$. For the first twist, the Koszul sequence
\[
0\to\OO_X\to\OO_X(A)\oplus\OO_X(\zeta)
\to\II_Z(\zeta+A)\to0
\]
and $H^1(S,\OO_S)=0$ give
$h^0(X,\II_Z(\zeta+A))=\nu_1+h^0(S,\VV)-1$.
For the twists by $-bA$ and $bA$, use
\[
\begin{aligned}
H^0(X,\II_Z(\zeta-2bA))&\subset H^0(S,\VV(-2bH))=0,\\
h^0(X,\II_Z(\zeta))&=1+h^0(S,\VV(-H)).
\end{aligned}
\]
The latter equality follows by twisting \eqref{eq:orig:4.3} by
$\zeta$ and using $H^i(S,\OO_S(-H))=0$ for $i=0,1$.
For the fourth twist, acyclicity in relative degree $-1$ gives
$H^0(X,F(-\zeta+(b+1)A))\simeq H^0(X,\II_Z(A))\simeq\C$.
Substituting yields
\[
\begin{aligned}
\dim\Ext^1_X(F,F)
&=\nu_{2b+1}-\nu_{2b}+\nu_1+h^0(S,\VV)\\
&\quad-h^0(S,\VV(-H))-3+\dim\ker\mu_F.
\end{aligned}
\]
Since $\nu_{2b+1}-\nu_{2b}=(4b+1)h$ and $\nu_1=h+2$, the
right-hand side equals
\[
(4b+2)h+h^0(S,\VV)-h^0(S,\VV(-H))-1+\dim\ker\mu_F
=D_b+\dim\ker\mu_F.
\]
This proves the formula and the pointwise equivalence.
Replacing a lift $\tau$ by $\tau+u\lambda$ precomposes $\mu_F$
with the invertible map $(\alpha,\beta)\mapsto
(\alpha,\beta+\lambda\alpha)$, so its injectivity is independent
of that choice.

For the component assertion, equality of local dimension and
tangent dimension will imply that the ambient local ring is regular.
The reduced closure $\mathcal M_b$ is irreducible of dimension $D_b$
by Lemma~\ref{orig:6.1}.
If $\mu_F$ is injective at a general member, upper semicontinuity
of $\dim\Ext^1(F,F)$ gives a dense open subset of $\mathcal T_b$ on which the
tangent dimension is $D_b$.
Its constructible image contains a dense open subset of $\mathcal M_b$.
At a closed point $y=[F]$ in that subset,
$\dim T_yM=D_b=\dim\mathcal M_b$. Lemma~\ref{lem:regularity}
gives the component and generic smoothness assertions.
Conversely, if $\mathcal M_b$ is a component and the ambient space is smooth
along its general points, its tangent dimension there is
$\dim \mathcal M_b=D_b$, and the formula forces $\ker\mu_F=0$.
The general singularity assertion follows from
Lemma~\ref{orig:4.4} and the nonempty open subset of reduced zero divisors.
\end{proof}

\subsection{Elementary transforms and the tangent correction}
For large $b$, the cokernel of $t\mapsto(s\wedge t)|_C$ measures
the excess tangent dimension. This separates the geometry of the
pair $(s,C)$ from the positivity needed for the Serre extension.

Hypothesis \eqref{eq:N} is not imposed in this subsection unless
explicitly stated.

Fix an admissible pair $(s,C)$: the curve $C=(u=0)\in|H|$ is
smooth and $s|_C$ is nowhere zero. Write $i:C\hookrightarrow S$.
\begin{definition}\label{orig:6.9}\label{def:elementary-transform}
The elementary transform associated with the pair $(s,C)$ is
\[
G=\ker\bigl(\VV\xrightarrow{s\wedge-}i_*((\det\VV)|_C)\bigr).
\]
\end{definition}
This is a variant of the Lazarsfeld--Mukai construction
\cite{Lazarsfeld}: the evaluation source $H^0(C,A_C)\otimes\OO_S$
for a globally generated line bundle $A_C$ on $C$ is replaced by
$\VV$, and the quotient on $C$ is $(\det\VV)|_C$.
The morphism is restriction to $C$ followed by exterior product
with $s$, and is surjective because $s$ is nowhere zero along $C$.
Completing $s$ to a local frame of $\VV$ identifies $G$ locally
along $C$ with $\OO_S\oplus\OO_S(-C)$, proving that $G$ is locally
free. Set
\[
\delta(s,C):=\dim\coker\!\bigl[H^0(S,\VV)\xrightarrow{s\wedge-}
H^0(C,(\det\VV)|_C)\bigr].
\]

Throughout this subsection, we say that $b$ is sufficiently large if
\[
b\ge1,\qquad 2bd>c_1(\VV)\cdot H,\qquad
H^1(S,G^\vee(2bH))=0.
\]
The numerical inequality also implies $H^0(S,\VV(-2bH))=0$,
since $\VV\hookrightarrow(\det\VV)^{\oplus h^0(S,\VV)}$.
The admissible-pair locus $U_0$ introduced in Section~\ref{subsec:extension-family} is defined
without \eqref{eq:N}. The corresponding universal quotient on
$\mathscr C_{U_0}\subset S\times U_0$, including the tautological
parameter twist, has a locally free kernel $\mathcal G$ with
fibers $G$. The bundle $\mathcal G^\vee$ is flat over the finite-type base
$U_0$. For $q:S\times U_0\to U_0$, relative Serre vanishing
\cite[III, Theorem~8.8]{HartshorneAG} yields a uniform bound $b_0$ such that
\[
R^iq_*\bigl(\mathcal G^\vee\otimes\pr_S^*\OO_S(2bH)\bigr)=0
\qquad(i=1,2,\ b\ge b_0).
\]
Here one may replace $H$ by a very ample multiple $rH$ and apply
the theorem to the finitely many twists by $jH$, $0\le j<r$.
Cohomology and base change \cite[III, Theorem~12.11]{HartshorneAG},
applied in descending degrees from the vanishing of fiber
cohomology above degree two, then gives
$H^1(S,G^\vee(2bH))=0$ for every pair in $U_0$ and $b\ge b_0$.
If $H$ is very ample and $G^\vee$ is $m$-regular, $2b\ge m-1$
suffices by Castelnuovo--Mumford regularity \cite[\S1.7]{HL}.
The correction $\delta(s,C)$ depends only on $\VV$ and $(s,C)$,
not on $b$. The required polarization inequalities are stated
separately in \textup{(iii)}.

\begin{proposition}\label{prop:5.7}
For all sufficiently large $b$, the following properties hold:
\begin{enumerate}[label=\textup{(\roman*)}]
\item For a fixed pair $(s,C)$ as above and every nonzero extension
$F=F_\sigma$, one has
\[
H^1(X,F(-bA))=0\qquad\text{and}\qquad
H^1(X,F(bA))\simeq H^1(S,G).
\]

\item If $F=F_\sigma$ is stable, then
\begin{align*}
\dim\Ext^1_X(F,F)
&=D_b+\dim\ker\!\bigl[H^1(S,G)\xrightarrow{\cdot u}H^1(S,G(H))\bigr]\\
&=D_b+\delta(s,C).
\end{align*}

\item Let $L\in\mathcal A_+$ be a rational polarization satisfying
$2bAL^2<\zeta L^2<2(b+1)AL^2$
\textup{(}or any rational $L\in\CC_b$ under~\eqref{eq:N}\textup{)}.
Then the image $\FF_b$ of the original Serre family has dimension $D_b$.
Its reduced closure $\mathcal M_b$ is an irreducible component
of $M_{X,L}(w)_{\mathrm{red}}$ along which $M_{X,L}(w)$ is generically
smooth if and only if $\delta(s,C)=0$ for a general pair $(s,C)$.
\end{enumerate}
\end{proposition}

\begin{proof}
To compute $\pi_*\II_Z(\zeta)$, twist the Koszul resolution
\eqref{eq:orig:4.3} by $\zeta$ to obtain
\[
0\to\OO_X(-A)\to\OO_X\oplus\OO_X(\zeta-A)
\to\II_Z(\zeta)\to0.
\]
The locally free terms have relative degrees $0,0,1$ and hence
vanishing higher direct images by Lemma~\ref{orig:2.8}.
\[
0\to\OO_S(-H)\to\OO_S\oplus\VV(-H)
\to\pi_*\II_Z(\zeta)\to0.
\]
The ideal sequence and $\OO_Z(\zeta)\simeq(\det\VV)|_C$ identify
\[
\pi_*\II_Z(\zeta)=\ker[\VV\twoheadrightarrow i_*((\det\VV)|_C)]=G,
\qquad R^j\pi_*\II_Z(\zeta)=0\quad(j>0).
\]
Pushing forward the three twists of the Serre sequence gives
\begin{align*}
0&\to\OO_S\to\pi_*F(-bA)\to G(-2bH)\to0,\\
0&\to\OO_S(2bH)\to\pi_*F(bA)\to G\to0,\\
0&\to\OO_S((2b+1)H)\to\pi_*F((b+1)A)\to G(H)\to0.
\end{align*}
Their higher direct images vanish. Kodaira vanishing gives
$H^i(\OO_S(nH))=0$ for $i=1,2$ and $n=2b,2b+1$.
The assumed vanishing and Serre duality give
\[
H^1(G(-2bH))\simeq H^1(G^\vee(2bH))^\vee=0.
\]
Together with $H^1(\OO_S)=0$, these prove \textup{(i)} and give
$H^1(F((b+1)A))\simeq H^1(G(H))$.

For \textup{(ii)}, apply $\Hom_X(-,F)$ to the locally free resolution \eqref{eq:orig:6.1}, with
$P=\OO_X(-(b+1)A)$ and
$Q=\OO_X(bA)\oplus\OO_X(-bA)\oplus
\OO_X(\zeta-(b+1)A)$. Its middle component $P\to\OO_X(-bA)$
is multiplication by $u$. The preceding vanishings and the
calculation for line bundles of relative degree $-1$ in Theorem~\ref{thm:regular-component} give
\[
\Ext^1(Q,F)=H^1(G),\qquad \Ext^1(P,F)=H^1(G(H)).
\]
Functoriality identifies the induced map with multiplication by $u$;
applying $\Hom(-,F)$ gives
\[
0\to\coker[\Hom(Q,F)\to\Hom(P,F)]\to\Ext^1(F,F)
\to\ker[H^1(G)\xrightarrow{u}H^1(G(H))]\to0.
\]
Our hypotheses give $H^0(\VV(-2bH))=0$ and the vanishing of higher cohomology for the relevant positive line bundles. The four $h^0$ counts in Theorem~\ref{thm:regular-component}, which do not
use $H^1(F(bA))=0$, give cokernel dimension $D_b$.
This gives the first equality in \textup{(ii)}.

We identify this kernel using restriction to $C=(u=0)$.
In the restriction of the defining sequence of $G$, the Tor term
is $(\det\VV)|_C\otimes\OO_C(-H)$, while the image of
$G|_C\to\VV|_C$ is the line subbundle $\OO_Cs$. Thus
\[
0\to(\det\VV)|_C\otimes\OO_C(-H)\to G|_C\to\OO_C\to0.
\]
The section $s$ splits this sequence, since its image in the last term
is $1$. Hence $G(H)|_C\simeq\OO_C(H)\oplus(\det\VV)|_C$.
Let $r_C$ be the restriction map
\[
r_C\colon H^0(S,G(H))\longrightarrow H^0(C,G(H)|_C),
\]
and denote the restriction--wedge map by
\[
\Phi_{s,C}\colon H^0(S,\VV)\longrightarrow H^0(C,(\det\VV)|_C),
\qquad t\longmapsto(s|_C)\wedge(t|_C).
\]
Choose Koszul signs so that the twisted presentation is
\[
0\to\OO_S\xrightarrow{(u,-s)}\OO_S(H)\oplus\VV
\xrightarrow{(s,u)}G(H)\to0.
\]
Since $H^1(S,\OO_S)=0$, every section of $G(H)$ can be written
as $as+ut$, with $a\in H^0(S,\OO_S(H))$ and $t\in H^0(S,\VV)$.
In a local frame $(s,e_2)$ of $\VV$, the sheaf $G$ has basis
$(s,ue_2)$. If $t=\alpha s+\beta e_2$, then $ut$ has coordinates
$(u\alpha,\beta)$ in this basis. Restricting to $C$ and identifying
$e_2$ with $s\wedge e_2$ in the determinant line therefore gives
\[
r_C(as+ut)=\bigl(a|_C,\Phi_{s,C}(t)\bigr)
\quad\text{in }H^0(C,\OO_C(H))\oplus H^0(C,(\det\VV)|_C).
\]
The map $H^0(S,\OO_S(H))\to H^0(C,\OO_C(H))$ is also surjective,
by $H^1(S,\OO_S)=0$. Consequently, under the displayed splitting,
\[
\operatorname{im}r_C=H^0(C,\OO_C(H))\oplus\operatorname{im}\Phi_{s,C}.
\]
The long exact sequence of $0\to G\xrightarrow{u}G(H)\to G(H)|_C\to0$
now yields
\[
\ker\bigl[H^1(S,G)\xrightarrow{u}H^1(S,G(H))\bigr]
\simeq\coker r_C\simeq\coker\Phi_{s,C},
\]
whose dimension is $\delta(s,C)$.
This proves \textup{(ii)} without assuming $H^1(G)=0$ or
$H^1(\VV)=0$.

For \textup{(iii)}, Theorem~\ref{thm:arbitrary-stability} gives stability for rational
$L\in\mathcal A_+$ satisfying $2bAL^2<\zeta L^2<2(b+1)AL^2$.
Lemma~\ref{lem:serre-family} applies to the Hilbert family used in
Lemma~\ref{orig:6.1}, giving a rank-$R_b$ extension bundle and an
injective classifying morphism.
The Hilbert-image dimension calculation is independent of \eqref{eq:N}, so
\[
\dim\mathcal M_b=D_b,
\qquad \dim T_{[F]}M=D_b+\delta(s,C).
\]
Lemma~\ref{lem:regularity}, together with the converse argument in the
proof of Theorem~\ref{thm:regular-component}, proves both implications
of the component criterion.
Under \eqref{eq:N}, these inequalities define $\CC_b$, so the conclusion holds for
every rational $L\in\CC_b$.
\end{proof}

Under the hypotheses of Proposition~\ref{prop:5.7}\textup{(ii)},
$H^1(S,G)=0$ is a sufficient condition for $\delta(s,C)=0$.

\subsection{Split \texorpdfstring{$\PP^1$}{P1}-bundles}\label{subsec:split}
For a split bundle, the elementary transform is a sum of line
bundles. The K3 line-bundle vanishings imply the component
criterion directly.

\par\medskip\noindent\refstepcounter{lemma}\label{orig:6.7}\textbf{Corollary~\ref{orig:6.7}.}
\textit{Assume \eqref{eq:N} and $\VV=\OO_S\oplus\OO_S(eH)$, with $e\ge0$.
For $b\ge1$, $2b>e$, and any rational $L\in\CC_b$, the reduced closure
of the Serre family is an irreducible component of $M_{X,L}(w)_{\mathrm{red}}$
along which $M_{X,L}(w)$ is generically smooth. Its dimension is
\[
D_b=
\begin{cases}
(4b+2)h+1,&e=0,\\
(4b+3)h+1,&e=1,\\
(4b+2e+1)h,&e\ge2.
\end{cases}
\]
A general member has reduced singularity scheme of length $(4b+2)h$.}
\begin{proof}
For a general section $s=(c,w)$, we have $c\ne0$. The automorphism
$(x,y)\mapsto(x/c,y-(w/c)x)$ of $\VV$ sends $s$ to $(1,0)$,
including when $e=0$. Hence
\[
G=\OO_S\oplus\OO_S((e-1)H),\qquad
R\pi_*\II_Z(\zeta)=G.
\]
The pushed-forward twists of the Serre sequence give
\[
\begin{aligned}
0&\to\OO_S\to\pi_*F(-bA)\to G(-2bH)\to0,\\
0&\to\OO_S(2bH)\to\pi_*F(bA)\to G\to0.
\end{aligned}
\]
All their higher direct images vanish. Since $H^1(\OO_S(nH))=0$
for every $n\in\Z$, these sequences imply
$H^1(F(-bA))=H^1(F(bA))=0$.
Theorem~\ref{thm:regular-component} now gives the component and
singularity assertions. The dimension formula follows from
Lemma~\ref{orig:6.1} and $h^0(\OO_S(nH))=2+n^2h$ for $n>0$.
\end{proof}

In the split case $\VV=\OO_S\oplus\OO_S(eH)$ with $e\ge0$,
assume \eqref{eq:N}, $b\ge1$, $2b>e$, and a rational polarization
$L\in\CC_b$. We retain $M=M_{X,L}(w)$ and $\mathcal M_b\subset M$
as defined at the start of Section~\ref{sec:components}.

\begin{proposition}\label{prop:6.3}
The original Serre family and its Brill--Noether loci satisfy the
following properties:
\begin{enumerate}[label=\textup{(\roman*)}]
\item Any member $F$ of the original Serre family whose
section $s$ has a nonzero constant component satisfies
\[\tag{5.5}\label{eq:6.1}
R^1\pi_*F=0,\qquad
\pi_*F\simeq\OO_S(bH)\oplus\OO_S(-bH)\oplus\OO_S((e-b-1)H).
\]
Consequently, $H^1(F)=H^3(F)=0$ and
\[\tag{5.6}\label{eq:6.2}
h^0(F)=k_b:=\nu_b+\nu_{e-b-1},\qquad
h^2(F)=\nu_b+\nu_{b+1-e}.
\]
In particular, if $e\le b$, then $k_b=b^2h+2$.

\item For every $1\le k\le k_b$, the variety
$\mathcal M_b$ is a $D_b$-dimensional irreducible
component of $(W^k)_{\mathrm{red}}$ along which $W^k$ is generically
smooth. Scheme-theoretically, $\mathcal M_b$ is entirely
contained in $W^k$ for all $1\le k\le k_b$, whereas the intersection
$W^{k_b+1}\cap\mathcal M_b$ is a (possibly empty)
proper closed subscheme of $\mathcal M_b$.
\end{enumerate}
\end{proposition}

\begin{proof}
For \textup{(i)}, a split-bundle automorphism sends a section $s$
with nonzero constant first component to $(1,0)$; when $e=0$, use
a constant change of basis. Its zero divisor on $X$ is the section
corresponding to $\VV\twoheadrightarrow\OO_S(eH)$.
Pushing forward the ideal sequence of $Z\simeq C$, twisted by
$\zeta$, gives
\[
\pi_*\II_Z(\zeta)\simeq\OO_S\oplus\OO_S((e-1)H),
\qquad R^1\pi_*\II_Z(\zeta)=0,
\]
since the map $\VV\to i_*\OO_C(eH)$ is projection followed by
restriction and is surjective. The pushed-forward extension splits
because
\[
\begin{aligned}
&\Ext^1_S\bigl(\OO_S(-bH)\oplus\OO_S((e-b-1)H),\OO_S(bH)\bigr)\\
&\qquad\simeq H^1(S,\OO_S(2bH))\oplus H^1(S,\OO_S((2b+1-e)H))=0.
\end{aligned}
\]
Since $H^1(S,\OO_S(aH))=0$ for every $a\in\Z$, the decomposition gives
$H^1(F)=H^3(F)=0$ and the stated value of $h^0(F)$.
Serre duality gives $h^2(\OO_S(nH))=\nu_{-n}$, so the same
decomposition yields $h^2(F)=\nu_b+\nu_{b+1-e}$, using $b\ge1$.
If $e\le b$, then
$h^0(S,\OO_S((e-b-1)H))=0$, so \eqref{eq:6.2} gives $k_b=b^2h+2$.

To prove \textup{(ii)}, let
$\mathcal J_k\subset\OO_M$
be the ideal sheaf of $W^k$, with intersections taken
scheme-theoretically in $M$. Upper semicontinuity and
\eqref{eq:6.2} give $h^0(F)\geq k_b$ for every $[F]\in \mathcal M_b$,
with equality on a dense open subset. Since $\mathcal M_b$ is reduced,
\[
\mathcal J_k\OO_{\mathcal M_b}=0\quad(1\leq k\leq k_b),
\qquad W^{k_b+1}\cap \mathcal M_b\subsetneq \mathcal M_b.
\]
By Corollary~\ref{orig:6.7}, $\mathcal M_b$ is an irreducible component
of $M_{\mathrm{red}}$ of dimension $D_b$; hence it is an irreducible
component of $(W^k)_{\mathrm{red}}$ for $1\leq k\leq k_b$.

At a general point $[F]\in\FF_b\subset\mathcal M_b$, the scheme $M$ is smooth,
$h^0(F)=k_b$, and $H^1(F)=0$. A minimal local perfect complex
representing the cohomology of the universal sheaf, on an
\etale\ chart, has degree-zero term $\OO^{k_b}$ and
zero degree-one term. The determinantal definition yields
$W^k=M$ locally at $[F]$ for $1\leq k\leq k_b$.
Thus $W^k$ is smooth at a general point of $\mathcal M_b$.
\end{proof}

\section{Relative Quot families of curve quotients}\label{sec:quot}
We enlarge the Serre family by allowing the kernel of a curve
quotient to vary among line bundles of degree zero. On a K3 surface, the condition
$e\ge2$ ensures the vanishing of the fixed-curve obstruction group; comparing
dimensions then identifies a generically smooth component.

Let $U_H=|H|^{\mathrm{sm}}$ be the open subset of smooth members of $|H|$,
and let $q\colon\mathscr C\to U_H$ be the universal curve.
Let $\mathscr Q\to U_H$ be the open locus in the relative Quot scheme
of $\VV$ on $\mathscr C$
parametrizing line-bundle quotients $\VV|_C\twoheadrightarrow J$
with kernel $K$ of degree zero. This degree fixes the Hilbert
polynomial of $J$, and local freeness of the quotient is an open
condition. For a quotient induced by an admissible pair $(s,C)$,
choose an irreducible component
$\mathscr Q_0\subset\mathscr Q$ containing it.

\begin{proposition}\label{orig:6.13}
Assume $\Pic(S)=\Z H$, $H^2=2h$, and
$\det\VV=\OO_S(eH)$ with $e\ge2$.
\begin{enumerate}[label=\textup{(\roman*)}]
\item The open locus $\mathscr Q$ is smooth of dimension $2eh+1$.
\item For all sufficiently large $b$ and all rational $L\in\CC_b$,
the component $\mathscr Q_0$, together with its Serre extensions,
gives a generically smooth stable moduli component of dimension
\[
D_b+\delta(s,C)=(4b+2e+1)h.
\]
\item Every member $F$ of this enlarged extension family satisfies
\[
h^0(F)=b^2h+2+h^0(S,\VV(-(b+1)H)).
\]
\end{enumerate}
\end{proposition}
\begin{proof}
For \textup{(i)}, the universal kernel on $\mathscr Q$ is a line
bundle of relative degree zero. For a fixed $C$, the relative tangent and
obstruction spaces at a quotient are
\[
\Hom_C(K,J)=H^0(C,J\otimes K^{-1}),\qquad
\Ext^1_C(K,J)=H^1(C,J\otimes K^{-1})
\]
\cite[Lemma~2.2.6 and Proposition~2.2.7]{HL}.
Since $J\otimes K^{-1}=(\det\VV)|_C\otimes K^{-2}$
has degree $2eh>2h$, Serre duality and Riemann--Roch give
\[
h^0(C,J\otimes K^{-1})=(2e-1)h,\qquad
H^1(C,J\otimes K^{-1})=0.
\]
Thus $\mathscr Q\to U_H$ is smooth of relative dimension $(2e-1)h$.
Since $U_H$ is smooth of dimension $h+1$, the irreducible components
of $\mathscr Q$ are disjoint and each has dimension
$(2e-1)h+(h+1)=2eh+1$.

For \textup{(ii)}, projectivizing the universal line-bundle quotient
$J$ identifies $\ZZ=\PP(J)\subset X\times\mathscr Q_0$ with
$\mathscr C\times_{U_H}\mathscr Q_0$, with $\OO_{\ZZ}(\zeta)\simeq J$.
Thus $\ZZ\to\mathscr Q_0$ is smooth and flat, with fibers $Z\simeq C$.
The embedded curve determines $C=\pi(Z)$ and the quotient
$\VV|_C\twoheadrightarrow\OO_Z(\zeta)$ by evaluation.
Lemma~\ref{lem:subsheaf-test}\textup{(a)} gives stability for every
nonzero extension. Thus Lemma~\ref{lem:serre-family} yields a smooth
irreducible parameter space $T$ with injective classifying morphism
$f:T\to M$, of dimension
$\dim\mathscr Q_0+R_b-1=(4b+2e+1)h$.

To identify its reflexive members, the local Koszul calculation gives
\[
\shExt^1_X(\II_Z(\zeta-bA),\OO_X(bA))
\simeq\det N_{Z/X}\otimes\OO_Z(2bA-\zeta)\simeq\mathcal L_K,
\]
where $Z\simeq C$. The normal sequence for $Z\subset X_C\subset X$
gives $\det N_{Z/X}\simeq\OO_Z(\zeta+A)\otimes K^{-1}$.
As in Lemma~\ref{orig:4.2}, the acyclicity of $\OO_X(2bA-\zeta)$
identifies global extensions with sections of this local extension
sheaf. The argument of Lemma~\ref{orig:4.4} therefore applies with
$\mathcal L_K$ in place of $\OO_C((2b+1)H)$.
Since $\deg\mathcal L_K=(4b+2)h\ge2g(C)+1$, it is very ample;
a general section has reduced zeros. Let $T^\circ\subset T$ be the
nonempty open subset of the corresponding reflexive sheaves with
reduced singularity scheme. It contains general members of the original
Serre family and has the same dimension as $T$.

To prove that the image closure is a generically smooth component,
we compare its dimension with the tangent dimension at a member
of the original Serre family. For a member associated with the chosen pair $(s,C)$,
Proposition~\ref{prop:5.7} gives tangent dimension $D_b+\delta$.
The kernel of the restriction--wedge map is
$H^0(S,G)\simeq\C\oplus H^0(S,\VV(-H))$, by the Koszul sequence and
$H^1(\OO_S(-H))=0$, so
$\delta=(2e-1)h-h^0(\VV)+1+h^0(\VV(-H))$.
Equivalently, $D_b+\delta=\dim\mathscr Q_0+R_b-1=(4b+2e+1)h$;
see Remark~\ref{rem:delta-codim}.
Let $\mathcal M_b'$ be the reduced closure of $f(T^\circ)$, and let $M$ be
the stable moduli space. The quotient induced by the chosen pair
lies in $\mathscr Q_0$, so a general extension associated with that pair
determines a point $y=[F]\in \mathcal M_b'$.
Quasi-finiteness gives $\dim \mathcal M_b'=\dim T^\circ=(4b+2e+1)h$.
At this point $\dim T_yM=D_b+\delta=\dim\mathcal M_b'$, so
Lemma~\ref{lem:regularity} shows that $\mathcal M_b'$ is an
irreducible component of $M_{\mathrm{red}}$ along which $M$ is
generically smooth. Finally, \textup{(iii)} is \eqref{eq:orig:6.8}, since
$h^0(\OO_S(bH))=b^2h+2$ on the K3 base.
\end{proof}

\begin{remark}\label{rem:delta-codim}
For $e\ge2$, the correction $\delta(s,C)$ has a direct interpretation
in the relative Quot scheme. The admissible-pair space $U_0$ of
Section~\ref{subsec:extension-family} maps to $\mathscr Q$ by taking
the quotient with kernel $\OO_Cs$. Its image lies in $\mathscr Q_0$:
$U_0$ is irreducible and the components of the smooth scheme
$\mathscr Q$ are disjoint. An embedded curve $Z\simeq C$ recovers
its quotient by evaluation, so this map and the Hilbert map $\eta_{\mathrm{Hilb}}$
have the same geometric fibers. Thus the section-induced quotient
locus has dimension $\dim B=h+h^0(\VV)-h^0(\VV(-H))$, and its
closure has codimension
\[
2eh+1-\bigl(h+h^0(\VV)-h^0(\VV(-H))\bigr)=\delta(s,C)
\]
in $\mathscr Q_0$. The excess tangent dimension in Theorem~\ref{thm:A}
therefore equals this codimension, while the enlarged Serre family of
Theorem~\ref{thm:B} has dimension $\dim\mathscr Q_0+R_b-1=D_b+\delta$.

For $\VV=\OO_S\oplus\OO_S(eH)$, every point of $\mathscr Q$ has
$K\simeq\OO_C$. Indeed, a nontrivial degree-zero line bundle $K$
has $\Hom_C(K,\OO_C)=0$, so its inclusion in $\VV|_C$ factors
through $\OO_C(eH)$. The resulting quotient has nonzero torsion
at the zeros of this map, since $\deg\OO_C(eH)=2eh>0$,
contrary to the definition of $\mathscr Q$. The torsion quotients
constructed in Section~\ref{subsec:second-component} therefore lie outside
$\mathscr Q$.
\end{remark}

The following example contrasts two behaviors for a nonsplit
bundle: the original Serre family may already fill a generically smooth
component, or it may have codimension one in the larger component
obtained from curve quotients.

\Needspace{4\baselineskip}
\begin{example}[Two nonsplit bundles on a double-plane K3 surface]\label{orig:6.14}
Let $f:S\to\PP^2$ be the double cover branched along a very general
smooth sextic, and put $H=f^*\OO(1)$. Then $\Pic(S)=\Z H$ by the
Noether--Lefschetz description for double-plane K3 surfaces
\cite[Chapter~17, \S1]{HuybrechtsK3}; explicit examples over
$\mathbf Q$ with geometric Picard rank one are given in \cite{EJ08}.
Thus $h=1$ and $f_*\OO_S=\OO\oplus\OO(-3)$. The projection formula
$f_*f^*E\simeq E\oplus E(-3)$ and the Euler sequence give
\begingroup
\renewcommand{\arraystretch}{1.25}
\setlength{\arraycolsep}{7pt}
\[
\begin{array}{c|ccccc}
\VV&e&\gamma&h^0(S,\VV)&h^0(S,\VV(-H))&\delta\\\hline
f^*T_{\PP^2}(-1)&1&2&3&0&0\\
f^*T_{\PP^2}&3&6&8&3&1
\end{array}
\]
\endgroup
Both bundles are globally generated, and their Chern classes exclude
a splitting into line bundles. For $C=f^{-1}(\lambda)$, where
$\lambda\subset\PP^2$ is a general line, the tangent--normal
sequence splits because $\Ext^1_\lambda(\OO_\lambda(1),\OO_\lambda(2))=0$.
Thus $T_{\PP^2}|_\lambda\simeq\OO_\lambda(2)\oplus\OO_\lambda(1)$,
and the two restrictions to $C$ are $\OO_C\oplus\OO_C(H)$ and
$\OO_C(H)\oplus\OO_C(2H)$, respectively. Restriction of global sections to
$C$ is surjective in both cases. The kernel of $s\wedge-$ is
$H^0(S,G)$, of dimension $1+h^0(S,\VV(-H))$.
For $\VV=f^*T_{\PP^2}(-1)$, the $e\ge2$ formula in
Proposition~\ref{orig:6.13} does not apply: the image has dimension
$3-1=2=h^0(\OO_C(H))$, so $\delta=0$ directly.
For $\VV=f^*T_{\PP^2}$, its dimension is
$8-1-3=4<5=h^0(\OO_C(3H))$, so $\delta=1$.
Since restriction is surjective, the sequence
$0\to\OO_C\xrightarrow{s}\VV|_C\to(\det\VV)|_C\to0$ also gives
\[
\delta(s,C)=\dim\ker\bigl[H^1(\OO_C)\xrightarrow{s}H^1(\VV|_C)\bigr].
\]
For the second bundle, this is the kernel of a surjection from a
two-dimensional space to a one-dimensional space; it describes the
infinitesimal variation of the degree-zero kernel $K$.

Fix $b$ sufficiently large and a rational $L\in\CC_b$. For the
first bundle, Proposition~\ref{prop:5.7} shows that the original
Serre family fills a generically smooth component of dimension
$4b+4$. For the second, the original family has dimension $4b+6$
but tangent dimension $4b+7$ at a general member, so its closure
is not a component along which the ambient moduli space is
generically smooth. Proposition~\ref{orig:6.13} enlarges it to
such a component of dimension $4b+7$; the extra direction varies
$K\in\Pic^0(C)$. In both cases \eqref{eq:orig:6.8} gives the
constant section count $b^2+2$ on the constructed families.
\end{example}

\begin{remark}\label{rem:double-plane-b}
In Example~\ref{orig:6.14}, the required large-$b$ vanishings
already hold for $b\ge2$. Indeed,
$0\to\OO_S\to G\to\II_{(s)_0}((e-1)H)\to0$, where $(s)_0$ is the
preimage of one point, respectively of three noncollinear points,
of $\PP^2$. For a general section, these points avoid the branch
sextic. Since $G^\vee\simeq G((1-e)H)$, twisting gives
\[
0\to\OO_S((2b+1-e)H)\to G^\vee(2bH)
\to\II_{(s)_0}(2bH)\to0.
\]
For $b\ge2$ and $e=1,3$, the first term has vanishing higher
cohomology. Thus $H^1(G^\vee(2bH))=0$ is equivalent to
surjectivity of evaluation of $H^0(\OO_S(2bH))$ on $(s)_0$.
The invariant and anti-invariant summands of $f_*\OO_S(2bH)$
separate values on the two sheets over these points. Evaluation
therefore reduces to the corresponding points of $\PP^2$ in
degrees $2b$ and $2b-3$, both surjective for $b\ge2$.
\end{remark}

\section{The boundary divisor}\label{sec:boundary}
Throughout this section, we assume \eqref{eq:N} and take
$\VV=\OO_S\oplus\OO_S(eH)$ with $e\ge2$. We fix $b\ge e+1$ and a
rational polarization $L\in\CC_b$; $M$ and
$\mathcal M_b$ are as in Section~\ref{subsec:split}.
The bound $b\ge e+1$ ensures both the inequality $b-e>0$, needed
for the Koszul restriction calculation in Lemma~\ref{lem:boundary-cohomology},
and the bound $2e-2b+1<0$ used in Lemma~\ref{lem:boundary-sections}.
We study the specialization $s_t=(t,v)$ to $s_0=(0,v)$.
Choose
$u\in H^0(S,\OO_S(H))$, $v\in H^0(S,\OO_S(eH))$, and
$\tau\in H^0(S,\OO_S((2b+1)H))$ such that $C=(u=0)$ is smooth
and $\Gamma=C\cap(v=0)$ and $\Xi=C\cap(\tau=0)$ are reduced and disjoint.
Recall that $\ell=\length\Gamma=2eh$, and that $k_b=b^2h+2$ and
$D_b=(4b+2e+1)h$ in this setting.
Set $s=(0,v)$ and define $F_0$ by
\[
0\to\OO_X(-(b+1)A)\xrightarrow{(\tau,u,s)}
\OO_X(bA)\oplus\OO_X(-bA)\oplus\OO_X(\zeta-(b+1)A)
\to F_0\to0.\tag{7.1}\label{eq:orig:9.3}
\]
The curve $Z=(u=s=0)$ is a torsion-quotient example of
Proposition~\ref{orig:6.12}. Indeed,
\[
0\to\OO_C\xrightarrow{(0,v|_C)}\OO_C\oplus\OO_C(eH)
\to J\to0,\qquad J\simeq\OO_C\oplus\OO_\Gamma(eH),
\]
and $Z=\PP_C(J)$ is the union of the section corresponding to
$\OO_C$ and the vertical fibers over $\Gamma$. Its kernel
is $\OO_C$, of degree zero, so Proposition~\ref{orig:6.12} gives
\[
h^0(F_0)=k_b,\qquad
\Ext^1_X(\II_Z(\zeta-bA),\OO_X(bA))
\simeq H^0(C,\OO_C((2b+1)H)),
\]
with extension dimension $R_b$ and vanishing $\Ext^2$.

The proof of Theorem~\ref{thm:C} proceeds as follows. Starting from
the presentation \eqref{eq:orig:9.3}, cohomology on $S$ yields the jumps
in $h^1$, $h^2$, and the tangent dimension, together with the
surjectivity of the section-lifting map
$\rho\colon\Ext_X^1(F_0,F_0)\to H^1(F_0(-bA))$
(Lemma~\ref{lem:boundary-cohomology} and Proposition~\ref{orig:9.6}).
Next, the universal family near $[F_0]$ is expressed as an extension
of a pullback from $S$ by a fixed line bundle
(Lemma~\ref{lem:boundary-descent}), and Bogomolov's inequality on $S$
forces the section counts to remain locally constant, placing the
local Brill--Noether ideals in the nilradical
(Lemma~\ref{lem:boundary-sections}). The map $\rho$ then controls their
linear terms and allows us to identify the two ideals
(Lemma~\ref{lem:boundary-ideals} and Proposition~\ref{cor:k3-boundary-nilpotents}),
while varying the parameters produces the boundary divisor $\Delta_b$
(Corollary~\ref{cor:boundary-divisor}). Finally,
Subsection~\ref{subsec:second-component} allows arbitrary reduced
divisors $\Gamma'\subset C$ of degree $\ell$ in the torsion-quotient
construction, producing a second component $\mathcal N_b$ through
$\Delta_b$ and determining the exact local dimension
$\dim\OO_{M,[F_0]}=D_b+h$ at a general closed point of $\Delta_b$
(Proposition~\ref{prop:torsion-family} and Corollary~\ref{cor:second-component}).

\subsection{The sheaf \texorpdfstring{$F_0$}{F0} and its cohomology}
\label{subsec:F0}
Write $P=\OO_X(-(b+1)A)$ and let $Q$ be the middle bundle in
\eqref{eq:orig:9.3}.

\begin{lemma}\label{lem:boundary-cohomology}
The sheaf $F_0$ defined by \eqref{eq:orig:9.3} satisfies
\[
h^0(F_0)=k_b,\qquad h^0(F_0(-bA))=1,\qquad
h^1(F_0)=\ell,\qquad h^3(F_0)=0,
\]
and
\[
\begin{aligned}
h^1(F_0(-bA))&=\ell-1,&\qquad h^1(F_0(bA))&=h+1,\\
h^1(F_0((b+1)A))&=1,& h^1(F_0(-\zeta+(b+1)A))&=0.
\end{aligned}
\]
Moreover, $h^0(F_0(-\zeta+(b+1)A))=1$ and
\[
\dim\Hom_X(P,F_0)-\dim\Hom_X(Q,F_0)+1=D_b.
\]
The map $\Ext^1_X(Q,F_0)\to\Ext^1_X(P,F_0)$ induced by
\eqref{eq:orig:9.3} is zero.
\end{lemma}
\begin{proof}
Since the constant component of $s=(0,v)$ is zero, the map into
$\VV(-(b+1)H)$ has zero component in $\OO_S(-(b+1)H)$.
Pushing \eqref{eq:orig:9.3} forward therefore splits off this summand:
$R^1\pi_*F_0=0$ and
$\pi_*F_0\simeq\OO_S(-(b+1)H)\oplus E$, where
\[
0\to\OO_S(-(b+1)H)\xrightarrow{(\tau,u,v)}
\OO_S(bH)\oplus\OO_S(-bH)\oplus\OO_S((e-b-1)H)
\to E\to0.\tag{7.2}\label{eq:orig:9.5}
\]
The entries $\tau,u,v$ have no common zero, so $E$ is locally free.
In particular, $H^3(F_0)=0$.

Put $\Sigma_m=H^0(S,\OO_S(mH))$, with
$\Sigma_m=0$ for $m<0$.
Since every line bundle $\OO_S(mH)$ has vanishing first cohomology,
\eqref{eq:orig:9.5} gives
\[
h^0(F_0(tA))=\nu_{t+b}+\nu_{t-b}+\nu_{t+e-b-1}.
\tag{7.3}\label{eq:boundary-h0-twists}
\]
Thus $h^0(F_0)=k_b$ and $h^0(F_0(-bA))=1$.
Applying Serre duality to the degree-two cohomology map in
\eqref{eq:orig:9.5}, after twisting by $tH$, also gives
\[
H^1(F_0(tA))^\vee\simeq
\Lambda_t:=\frac{\Sigma_{b+1-t}}
{\tau \Sigma_{-b-t}+u\Sigma_{b-t}+v\Sigma_{b+1-e-t}}.
\tag{7.4}\label{eq:boundary-h1-cokernel}
\]

For $m=b+1,2b+1$, the inequality $m-e-1>0$ gives
$H^2(S,\OO_S((m-e-1)H))=0$. Together with the line-bundle
$H^1$ vanishings, the Koszul resolution of $\Gamma=(u=v=0)$ yields
\[
\frac{\Sigma_m}{u\Sigma_{m-1}+v\Sigma_{m-e}}
\simeq H^0(\Gamma,\OO_\Gamma(mH)).
\]
Consequently, $\dim \Lambda_0=\ell$.
Since $\tau|_\Gamma$ is nowhere zero, its class spans a one-dimensional
subspace, and hence $\dim \Lambda_{-b}=\ell-1$.
The remaining cases of \eqref{eq:boundary-h1-cokernel} are
\[
\Lambda_b=\Sigma_1/\C u,\qquad \Lambda_{b+1}=\C,
\]
of dimensions $h+1$ and $1$, respectively.

Twisting \eqref{eq:orig:9.3} by $-\zeta+(b+1)A$ leaves one
summand $\OO_X$; all other terms have relative degree $-1$.
Lemma~\ref{orig:2.8} gives
\[
R\pi_*F_0(-\zeta+(b+1)A)\simeq\OO_S.
\]
Hence $h^1(F_0(-\zeta+(b+1)A))=0$ and
$h^0(F_0(-\zeta+(b+1)A))=1$.
Using \eqref{eq:boundary-h0-twists}, we obtain
\[
\begin{aligned}
&\dim\Hom_X(P,F_0)-\dim\Hom_X(Q,F_0)+1\\
&\qquad=\nu_{2b+1}-\nu_{2b}+\nu_1+\nu_e-\nu_{e-1}-2=D_b.
\end{aligned}
\]

It remains to identify the map on $\Ext^1$.
The vanishing just proved reduces it to
\[
H^1(F_0(-bA))\oplus H^1(F_0(bA))
\xrightarrow{\,\tau+u\,}H^1(F_0((b+1)A)).
\]
For convenience, set
\[
\begin{aligned}
\Lambda_\circ&:=\Lambda_{b+1}=\C,\\
\Lambda_-&:=\Lambda_{-b}=\frac{\Sigma_{2b+1}}{\tau\C+u\Sigma_{2b}+v\Sigma_{2b+1-e}},\\
\Lambda_+&:=\Lambda_b=\Sigma_1/\C u.
\end{aligned}
\]
Naturality of Serre duality in \eqref{eq:boundary-h1-cokernel}
identifies the dual map with
$\Lambda_\circ\longrightarrow\Lambda_-\oplus\Lambda_+$,
$1\longmapsto([\tau],[u])=(0,0)$, since both classes vanish by the
definitions of the quotients. Thus the original map is zero.
\end{proof}

\Needspace{4\baselineskip}
\begin{proposition}\label{orig:9.6}
The sheaf $F_0$ defined by \eqref{eq:orig:9.3} is a stable reflexive
boundary point of $\mathcal M_b$.
Its reduced singularity scheme has length $(4b+2)h$, and
\[
h^0(F_0)=k_b,\qquad \dim\Ext^1_X(F_0,F_0)=D_b+(2e+1)h.
\]
Moreover,
\[
h^1(F_0)=\ell,\qquad
h^2(F_0)=4+(e^2-2eb+2b^2+2b+1)h,
\]
so $h^2$ increases by $\ell$ relative to the original Serre family.
The section-lifting map
$\rho\colon\Ext^1_X(F_0,F_0)\to H^1(F_0(-bA))\simeq\C^{\ell-1}$
is surjective.
\end{proposition}
\begin{proof}
Put $Z=(u=s=0)$. Since $u,s$ form a regular sequence,
\eqref{eq:orig:9.3} defines an extension
\[
0\to\OO_X(bA)\to F_0\to\II_Z(\zeta-bA)\to0.
\]
The curve $Z$ consists of the section corresponding to $\OO_C$
and the vertical fibers over $\Gamma$.
Along these fibers, $\tau|_\Gamma$ is a unit, so $F_0$ is locally free.
Its remaining nonfree points lie over $\Xi$, where the entries
of the local presentation form a regular system of parameters.
The local calculation of Lemma~\ref{orig:4.3} gives
\[
F_0\text{ reflexive},\qquad \Sing(F_0)\simeq\Xi,\qquad
\length\Sing(F_0)=(4b+2)h.
\]
The extension is nonzero, since its local class is represented
by $\tau|_Z$. Its Chern character is determined by
\eqref{eq:orig:9.3}. Every $L\in\CC_b$ satisfies
\eqref{eq:slope-window} by \eqref{eq:orig:5.1}, so
Lemma~\ref{lem:subsheaf-test}\textup{(a)} gives stability.

Replacing $s=(0,v)$ by $s_t=(t,v)$ gives a fiberwise injective map
between vector bundles over $X\times\mathbf A^1$.
Its cokernel is flat over $\mathbf A^1$ by the local criterion of
flatness. For $t\ne0$, the section $s_t$ is nowhere zero on $S$,
so $([s_t],C)\in U_0$ and the fiber belongs to the original Serre
family. Hence
\[
[F_0]\in(\overline{\FF_b})_{\mathrm{red}}=\mathcal M_b.
\]
Moreover, $h^0(F_0(-bA))=1$ by Lemma~\ref{lem:boundary-cohomology}.
Its distinguished line subsheaf is unique up to scalar, and its
quotient curve, which contains vertical fibers, distinguishes
$F_0$ from members of the original Serre family.

The flat family $s_t=(t,v)$ has constant Euler characteristic.
By Lemma~\ref{lem:boundary-cohomology} and
Proposition~\ref{prop:6.3}\textup{(i)}, $F_0$ and a general member
$F$ of the original family have the same $h^0$ and vanishing $h^3$,
while $h^1(F)=0$ and $h^1(F_0)=\ell$. Hence
\[
\begin{aligned}
h^2(F_0)&=h^2(F)+\ell=\nu_b+\nu_{b+1-e}+\ell\\
&=4+(e^2-2eb+2b^2+2b+1)h.
\end{aligned}
\]

Finally, apply $\Hom_X(-,F_0)$ to $0\to P\to Q\to F_0\to0$.
Stability gives $\Hom_X(F_0,F_0)=\C\cdot\mathrm{id}_{F_0}$.
Lemma~\ref{lem:boundary-cohomology} shows that the Hom cokernel
has dimension $D_b$ and that the subsequent map on $\Ext^1$ is zero.
Hence there is an exact sequence
\[
0\to\coker\bigl[\Hom_X(Q,F_0)\to\Hom_X(P,F_0)\bigr]
\to\Ext^1_X(F_0,F_0)\to\Ext^1_X(Q,F_0)\to0.
\]
Here $\Ext^1_X(Q,F_0)=H^1(F_0(-bA))\oplus H^1(F_0(bA))$, by
Lemma~\ref{lem:boundary-cohomology}. It follows that
\[
\dim\Ext^1_X(F_0,F_0)
=D_b+(\ell-1)+(h+1)=D_b+(2e+1)h.
\]
Composing the right-hand surjection with the first projection is
pullback along $\OO_X(bA)\hookrightarrow F_0$.
This is the section-lifting map, so it is surjective.
\end{proof}

The surjection $\rho$ of Proposition~\ref{orig:9.6} detects first-order
obstructions to lifting $\OO_X(bA)\hookrightarrow F_0$.
To identify the corresponding BN ideal, put
\[
\mathcal B_b=\{[F]\in M:h^0(F(-bA))\ge1\},
\]
and equip $\mathcal B_b$ with its determinantal structure.

\subsection{Local structure near \texorpdfstring{$[F_0]$}{[F0]}}
\label{subsec:local}
The following lemma expresses the quotient of the universal
family by $\OO_X(\zeta-(b+1)A)$ as a pullback from $S$, so that
section counts can be studied through extensions on the K3 surface.

Choose an \etale\ universal
chart $(U,t_0)\to(M,[F_0])$ with universal sheaf $\widetilde{\FF}$
and fibers $F_t$. Write $\pi_U:X\times U\to S\times U$.

\Needspace{4\baselineskip}
\begin{lemma}\label{lem:boundary-descent}
After shrinking $U$, the following hold.
\begin{enumerate}[label=\textup{(\roman*)}]
\item The universal sheaf fits into a $U$-flat exact sequence
\[
0\to\pr_X^*\OO_X(\zeta-(b+1)A)\to\widetilde{\FF}\to\pi_U^*\widetilde{\mathcal A}\to0,
\]
where $\widetilde{\mathcal A}$ is $U$-flat with geometric fibers
$\widetilde{\mathcal A}_t\simeq\II_{\Xi_t}((b+1)H)$ for
zero-dimensional subschemes $\Xi_t\subset S$ of length $(2b+1)d$.
\item For every geometric point $t\to U$, pushing forward gives
\[
0\to\OO_S(-(b+1)H)\oplus\OO_S((e-b-1)H)
\to\pi_*F_t\to\II_{\Xi_t}((b+1)H)\to0.
\]
Pushout along the projection onto $\OO_S((e-b-1)H)$ defines a
$U$-flat family with locally free fibers $E_t$ and exact sequences
\[
\begin{aligned}
0&\to\OO_S((e-b-1)H)\to E_t\to\II_{\Xi_t}((b+1)H)\to0,\\
0&\to\OO_S(-(b+1)H)\to\pi_*F_t\to E_t\to0.
\end{aligned}
\]
\item Every geometric fiber satisfies
$h^0(X,F_t(-bA))\le1$ and $h^0(X,F_t)\le k_b$.
\end{enumerate}
\end{lemma}
\begin{proof}
To construct a quotient that descends to $S\times U$, we first extend
$\OO_X(\zeta-(b+1)A)\hookrightarrow F_0$ to the universal family.
The quotient by this summand is the pullback of the cokernel of
$(\tau,u)\colon\OO_S(-(b+1)H)\to\OO_S(bH)\oplus\OO_S(-bH)$.
The Koszul complex of $(\tau,u)$ identifies that cokernel with
$\II_\Xi((b+1)H)$. Thus
\[
0\to \OO_X(\zeta-(b+1)A)\to F_0\to\pi^*\II_\Xi((b+1)H)\to0,
\]
and Lemma~\ref{lem:boundary-cohomology} gives
\[
\Hom_X(\OO_X(\zeta-(b+1)A),F_0)\simeq\C,\qquad \Ext^1_X(\OO_X(\zeta-(b+1)A),F_0)=0.
\]
A minimal representative of the relative Hom complex has ranks $1$
and $0$ in degrees zero and one. After shrinking $U$, it has
$K^0\simeq\OO_U$ and $K^1=0$, so its degree-zero cohomology is a
trivial line bundle and commutes with arbitrary base change. A trivialization defines an evaluation
map $\pr_X^*\OO_X(\zeta-(b+1)A)\to\widetilde{\FF}$ that is injective on the fiber over $t_0$.
The fiberwise flatness criterion and properness of $X$ then give,
after further shrinking, a $U$-flat exact sequence
\[
0\to\pr_X^*\OO_X(\zeta-(b+1)A)\to\widetilde{\FF}\to\widetilde{\QQ}\to0.
\]

The quotient $\widetilde{\QQ}$ is $U$-flat, with central fiber
$\pi^*\II_\Xi((b+1)H)$. Since $\OO_S\to R\pi_*\OO_X$ is an
isomorphism, Lemma~\ref{lem:pullback-descent} gives, after shrinking
$U$, a coherent $U$-flat sheaf
$\widetilde{\mathcal A}=\pi_{U*}\widetilde{\QQ}$ with
$\widetilde{\QQ}\simeq\pi_U^*\widetilde{\mathcal A}$, compatibly
with arbitrary base change over $U$.

After shrinking $U$ further, every geometric fiber
$\widetilde{\mathcal A}_t$ has rank one and is torsion-free on $S$,
by openness of purity \cite[Proposition~2.3.1]{HL}.
Its determinant is constant near $t_0$, because the Picard scheme
of $S$ is discrete and reduced by $H^1(S,\OO_S)=0$.
The fixed determinant $\OO_S((b+1)H)$ and constancy of the Hilbert
polynomial give
\[
\widetilde{\mathcal A}_t\simeq\II_{\Xi_t}((b+1)H),\qquad
\length\Xi_t=(2b+1)d.
\]
These identifications are fiberwise and require no base-change
statement for relative double duals.

This proves \textup{(i)}. Since $R^1\pi_*\OO_X(\zeta-(b+1)A)=0$, pushing forward
the sequence in \textup{(i)} gives the extension in \textup{(ii)}.
The construction commutes with base change. The pushout family
defining $E_t$ is an extension of $\widetilde{\mathcal A}$ by
$\pr_S^*\OO_S((e-b-1)H)$, hence is $U$-flat. At $t_0$, the sheaf
$E_{t_0}$ is the locally free bundle $E$ in
\eqref{eq:orig:9.5}. Openness of local freeness and properness of
$S$ give locally free $E_t$ after shrinking $U$.
Finally, upper semicontinuity and
Lemma~\ref{lem:boundary-cohomology} give \textup{(iii)}.
\end{proof}

\Needspace{4\baselineskip}
To turn the upper bounds in Lemma~\ref{lem:boundary-descent}\textup{(iii)}
into equalities, it suffices to construct a nonzero section of
$F_t(-bA)$. The Chern classes of $E_t$ force a nonzero section
of $E_t(-bH)$, which lifts to $F_t(-bA)$ by cohomology vanishing
on $S$.

\begin{lemma}\label{lem:boundary-sections}
Let $U\to M$ and $\widetilde{\FF}$ be the \etale\ morphism and
universal family constructed in Lemma~\ref{lem:boundary-descent}.
For every geometric point $t\to U$, the corresponding sheaf $F_t$
satisfies
\[
h^0(F_t(-bA))=1,\qquad h^0(F_t)=k_b.
\]
\end{lemma}
\begin{proof}
Fix a geometric point $t$ and write $F=F_t$.
By Lemma~\ref{lem:boundary-descent}\textup{(ii)}, $E_t$ is locally
free and has Chern classes
\[
c_1(E_t)=eH,\qquad c_2(E_t)=\bigl((b+1)e-b^2\bigr)d<0,
\]
where the inequality follows from $b\ge e+1$. Thus
$4c_2(E_t)-c_1(E_t)^2<0$, and Bogomolov's inequality
\cite[Theorem~3.4.1]{HL} implies that $E_t$ is not slope-semistable.
A saturated destabilizing subsheaf is reflexive, since $E_t$ is
locally free and its quotient is torsion-free. It is therefore a line
bundle, which we write as $\OO_S(aH)$ using $\Pic(S)=\Z H$.
The slope inequality gives $a>e/2$, and the quotient has the form
$\II_Y((e-a)H)$ for a zero-dimensional subscheme $Y\subset S$.
Comparing second Chern classes yields
\[
0\le\length Y=\bigl((b+1)e-b^2+a(a-e)\bigr)d.
\]
If $a\le b-1$, the function $a(a-e)$ is strictly increasing for
$a>e/2$, so
\[
\length Y\le\bigl((b+1)e-b^2+(b-1)(b-1-e)\bigr)d
=(2e-2b+1)d<0,
\]
a contradiction. Hence $a\ge b$, and the inclusion
$\OO_S(aH)\hookrightarrow E_t$ gives
\[
h^0(E_t(-bH))\ge h^0(\OO_S((a-b)H))\ge1.
\]

To lift such a section to $F(-bA)$, we twist the second exact sequence
of Lemma~\ref{lem:boundary-descent}\textup{(ii)} by $-bH$:
\[
0\to\OO_S(-(2b+1)H)\to(\pi_*F)(-bH)\to E_t(-bH)\to0.
\]
Since $H^i(S,\OO_S(-(2b+1)H))=0$ for $i=0,1$, the projection
formula and this sequence give
\[
H^0(X,F(-bA))\simeq H^0(S,E_t(-bH)).
\]
Thus $h^0(F(-bA))\ge1$. Together with the upper bound in
Lemma~\ref{lem:boundary-descent}\textup{(iii)}, this gives $h^0(F(-bA))=1$.
Its nonzero section injects $\OO_X(bA)$ into the torsion-free sheaf
$F$, so $h^0(F)\ge h^0(\OO_X(bA))=k_b$.
The other upper bound in that lemma then gives $h^0(F)=k_b$.
\end{proof}

\Needspace{4\baselineskip}
Lemma~\ref{lem:boundary-sections} implies that the ideals of
$\mathcal B_b\times_M U$ and $W^{k_b}\times_M U$ lie in the
nilradical of $\OO_U$. To identify their scheme structures, we
now compare the ideals themselves, using $\rho$ to determine
their linear terms.

\begin{lemma}\label{lem:boundary-ideals}
For the \etale\ morphism $U\to M$ of
Lemma~\ref{lem:boundary-descent}, put $(R,\mathfrak m)=\OO_{U,t_0}$.
Let $I,J\subset R$ be the ideals of the inverse images of
$\mathcal B_b$ and $W^{k_b}$, respectively. Then
\[
I=J\subset\sqrt{0},\qquad
T_{t_0}\mathcal B_b=T_{t_0}W^{k_b}=\ker\rho.
\]
The ideal $I$ has $\ell-1$ generators with linearly independent
images in $\mathfrak m/\mathfrak m^2$.
\end{lemma}
\begin{proof}
Lemma~\ref{lem:boundary-sections} implies
$I,J\subset\sqrt{0}$. The map
$\rho:\Ext^1(F_0,F_0)\twoheadrightarrow H^1(F_0(-bA))$
is surjective by Proposition~\ref{orig:9.6}.

We compute both tangent spaces using first-order deformations.
Write $\C[\epsilon]=\C\oplus\epsilon\C$, with $\epsilon^2=0$, for the
algebra of dual numbers. Let $\iota_0\in H^0(F_0(-bA))$ be the
section defined by the inclusion $\OO_X(bA)\hookrightarrow F_0$.
The \etale\ universal chart identifies
$T_{t_0}U$ with $\Ext^1_X(F_0,F_0)$: a tangent vector
$\xi\in(\mathfrak m/\mathfrak m^2)^\vee$ defines
$\varphi_\xi:R\to \C[\epsilon]$, with
$\varphi_\xi(f)=\epsilon\xi(\bar f)$ for $f\in\mathfrak m$.
Restricting the universal family along $\varphi_\xi$ gives a flat
first-order deformation $F_\xi$ whose underlying extension on $X$ is
\[
0\to F_0\xrightarrow{\ \epsilon\ }F_\xi\to F_0\to0
\]
and has class $\xi$ under this identification.

A minimal complex for the twisted universal family begins
\[
R\xrightarrow{\ d^0\ }R^{\ell-1}\longrightarrow\cdots,
\qquad d^0=(f_1,\ldots,f_{\ell-1})^{\mathsf t},\quad f_i\in\mathfrak m.
\]
Its degree-zero basis reduces to $\iota_0$, and the determinantal
definition gives $I=(f_1,\ldots,f_{\ell-1})$.
Derived base change identifies $H^0(F_\xi(-bA))$ with the kernel
of $d^0\otimes_{\varphi_\xi}\C[\epsilon]$. Reduction to
$H^0(F_0(-bA))=\C\iota_0$ sends $c_0+\epsilon c_1$ to $c_0\iota_0$, while
\[
(d^0\otimes_{\varphi_\xi}\C[\epsilon])(c_0+\epsilon c_1)
=\epsilon c_0\bigl(\xi(\bar f_1),\ldots,
\xi(\bar f_{\ell-1})\bigr)^{\mathsf t}.
\]
Thus $\iota_0$ lifts if and only if $\xi(\bar f_i)=0$ for every $i$,
equivalently $\xi\in T_{t_0}\mathcal B_b$.
Applying $\Hom_X(\OO_X(bA),-)$ to the extension above identifies
the same obstruction with the connecting image
$\iota_0^*\xi=\rho(\xi)$.
\[
T_{t_0}\mathcal B_b=\ker\rho.
\]
This tangent space is the common annihilator of the classes
$\bar f_i\in\mathfrak m/\mathfrak m^2$. Surjectivity of $\rho$
makes its codimension $\ell-1$, so these $\ell-1$ classes are
linearly independent.

The same argument applies to the untwisted family. Its minimal
complex has degree-zero term $R^{k_b}$, with a basis reducing to
$a_1\iota_0,\ldots,a_{k_b}\iota_0$, where $a_1,\ldots,a_{k_b}$ is a basis
of $H^0(\OO_S(bH))$. Thus $\xi\in T_{t_0}W^{k_b}$ if and only
if every $a_i\iota_0$ lifts to $F_\xi$. By naturality of the connecting
map, this is equivalent to $a_i\rho(\xi)=0$ in $H^1(F_0)$ for
all $i$. Choose $a\in H^0(\OO_S(bH))$ nonzero at every point
of $\Gamma$. By \eqref{eq:boundary-h1-cokernel} and the restriction
identifications in Lemma~\ref{lem:boundary-cohomology}, the dual
of multiplication by $a$ on $H^1$ is
\[
H^0(\Gamma,\OO_\Gamma((b+1)H))
\xrightarrow{\ \cdot a\ }
H^0(\Gamma,\OO_\Gamma((2b+1)H))/\C\tau.
\]
It is surjective since $a$ is nowhere zero on $\Gamma$.
Multiplication by $a$ is injective on
$H^1(F_0(-bA))$. Since $a$ is a linear combination of the $a_i$,
the vanishing of all $a_i\rho(\xi)$ forces $\rho(\xi)=0$, giving
\[
T_{t_0}W^{k_b}=\ker\rho=T_{t_0}\mathcal B_b.
\]

The tangent-space equality determines only the linear parts of
the ideals. To compare the ideals themselves, we use minimal
cohomology complexes, whose first differentials generate the two
BN ideals. Over $R/I$, derived base change identifies the
degree-zero generator with a section $\iota_I$ of the twisted family
lifting $\iota_0$.

Choose a basis $a_1,\ldots,a_{k_b}$ of $H^0(S,\OO_S(bH))$.
The cycles $a_i \iota_I$ reduce to a basis of $H^0(X,F_0)$.
Let $\tilde d^{\,0}$ be the first differential of a minimal complex for the
untwisted family, with degree-zero term $R^{k_b}$, and let $A_I$
be the coordinate matrix of these cycles over $R/I$. Since the
complex starts in degree zero, its zeroth cohomology after any
base change is the kernel of the first differential. Hence
\[
A_I\in\operatorname{GL}_{k_b}(R/I),\qquad
(\tilde d^{\,0}\bmod I)A_I=0
\quad\Longrightarrow\quad \tilde d^{\,0}\bmod I=0.
\]
By the determinantal definition, $J=I(W^{k_b})$ is generated by
the entries of $\tilde d^{\,0}$. Hence $J\subset I$ scheme-theoretically.

The tangent-space equality gives
$(I+\mathfrak m^2)/\mathfrak m^2=(J+\mathfrak m^2)/\mathfrak m^2$.
The independent linear parts of the generators $f_i$ imply
$I\cap\mathfrak m^2=\mathfrak m I$: if
$\sum r_i f_i\in\mathfrak m^2$, then every $r_i\in\mathfrak m$.
Together with $J\subset I$, this gives
\[
I=J+(I\cap\mathfrak m^2)=J+\mathfrak m I.
\]
Nakayama's lemma applied to $I/J$ yields $I=J$.
\end{proof}

\Needspace{4\baselineskip}
\begin{proposition}\label{cor:k3-boundary-nilpotents}
For the sheaf $F_0$ of Proposition~\ref{orig:9.6}, let
$\rho:\Ext^1_X(F_0,F_0)\to\C^{\ell-1}$ be the section-lifting map.
\begin{enumerate}[label=\textup{(\roman*)}]
\item There is an open neighborhood $\Omega\subset M$ of $[F_0]$ such that
\[
W^{k_b+1}\cap\Omega=\varnothing,\qquad
\Omega_{\mathrm{red}}=(W^{k_b}\cap\Omega)_{\mathrm{red}}
=(\mathcal B_b\cap\Omega)_{\mathrm{red}}.
\]
\item In $(R,\mathfrak m)=\OO_{M,[F_0]}$, the local ideals of
$W^{k_b}$ and $\mathcal B_b$ coincide in an ideal
$I\subset\sqrt{0}$ generated by $\ell-1=2eh-1$ elements
with linearly independent images in $\mathfrak m/\mathfrak m^2$.
Consequently, $M$ is nonreduced at $[F_0]$, every tangent vector
obtained from a reduced base lies in
\[
T_{[F_0]}W^{k_b}=T_{[F_0]}\mathcal B_b=\ker\rho,
\]
and
\[
\dim R\le\dim\ker\rho=D_b+h+1<\dim T_{[F_0]}M.
\]
\end{enumerate}
\end{proposition}
\begin{proof}
Let $\Omega\subset M$ be the open image of the \etale\ morphism
$U\to M$ constructed in Lemma~\ref{lem:boundary-descent}.
Assertion \textup{(i)} follows from Lemma~\ref{lem:boundary-sections}.

Put $R=\OO_{M,[F_0]}$, with maximal ideal $\mathfrak m$.
By Lemma~\ref{lem:boundary-ideals} and faithful flatness, the
local ideals of $W^{k_b}$ and $\mathcal B_b$ coincide:
\[
I(W^{k_b})=I(\mathcal B_b)=I\subset\sqrt{0}.
\]
Their common ideal has $\ell-1$ generators with linearly
independent images in $\mathfrak m/\mathfrak m^2$.
Since $\ell-1>0$, this ideal is nonzero, and $R$ is nonreduced.
The \etale\ identification of tangent spaces also gives
\[
T_{[F_0]}W^{k_b}=T_{[F_0]}\mathcal B_b=\ker\rho.
\]

The surjectivity of $\rho$ and Proposition~\ref{orig:9.6}
now imply
\[
\begin{aligned}
\dim\ker\rho
&=\dim\Ext^1_X(F_0,F_0)-(\ell-1)\\
&=D_b+(2e+1)h-(\ell-1)=D_b+h+1.
\end{aligned}
\]
Every morphism from a reduced scheme to $\Omega$ annihilates the
nilpotent ideal of $\mathcal B_b\cap \Omega$; its tangent image at
$[F_0]$ lies in $\ker\rho$.
Applying this to $\Omega_{\mathrm{red}}\hookrightarrow \Omega$ yields
\[
\dim R=\dim R_{\mathrm{red}}
\le\dim T_{[F_0]}\Omega_{\mathrm{red}}
\le\dim\ker\rho<\dim T_{[F_0]}M.
\]
\end{proof}

\subsection{The divisor \texorpdfstring{$\Delta_b$}{Delta b}}
\label{subsec:divisor}
The preceding arguments apply to every choice of $(u,v,\tau)$
with the stated smoothness and disjointness properties. To determine
the dimension of this boundary locus, we parametrize $v|_C$ up to
scalar, thereby removing the redundancy in its lift to $S$.

\begin{corollary}\label{cor:boundary-divisor}
The closure of the locus of sheaves $F_0$ in
Proposition~\ref{orig:9.6} is an irreducible divisor
$\Delta_b\subset\mathcal M_b$ of dimension $D_b-1$.
Its general point lies outside the original Serre family and satisfies
all the cohomological and local BN conclusions of
Propositions~\ref{orig:9.6} and~\ref{cor:k3-boundary-nilpotents}.
In particular, the ambient moduli scheme is nonreduced at a general
point of $\Delta_b$.
\end{corollary}
\begin{proof}
We construct a flat parameter family of dimension $D_b-1$ and
recover its parameters from the middle term.

For the universal curve $q\colon\mathscr C\to U_H$ of Section~\ref{sec:quot},
the sheaf $q_*\OO_{\mathscr C}(eH)$ is locally free of rank $(2e-1)h$
and commutes with base change, since $e\ge2$. Let $B_\partial$ be the
nonempty open subset of $\PP_{U_H}\bigl(q_*\OO_{\mathscr C}(eH)\bigr)$
parametrizing pairs $(C,[v_C])$ for which $\Gamma=(v_C)_0$ is reduced. Then
\[
\dim B_\partial=(h+1)+((2e-1)h-1)=2eh.
\]
Over $B_\partial$, the tautological image of $(0,v_C)$ cuts out
$Z=\PP_C(\OO_C\oplus\OO_\Gamma(eH))$. It is a nonzerodivisor on
each fiber of the universal surface, so these divisors form a flat
family of relative complete intersections of codimension two in
$X\times B_\partial$. The embedded curve recovers $C$ by projection
and $[v_C]$ from the divisor of its vertical fibers.
Lemma~\ref{lem:subsheaf-test}\textup{(a)} gives stability of all
nonzero extensions. We may therefore apply
Lemma~\ref{lem:serre-family}. In the projective extension bundle,
let $T_\partial$ be the nonempty open subset where
$\Xi=(\sigma)_0$ is reduced and disjoint from $\Gamma$.
The extension line bundle is $\OO_C((2b+1)H)$, so this open subset
is nonempty over every point of $B_\partial$.
It carries a flat family with injective classifying morphism and has
dimension $2eh+R_b-1=D_b-1$.

To identify its image, the restriction maps from
$H^0(S,\OO_S(eH))$ and $H^0(S,\OO_S((2b+1)H))$ to $C$ are
surjective. Every fiber therefore admits lifts $v,\tau$ and is a
sheaf $F_0$ of Proposition~\ref{orig:9.6}. Replacing $v$ by $v+uw$
or $\tau$ by $\tau+u\lambda$ changes the presentation by an
automorphism of $Q$; rescaling the lifts does so as well.
Thus $F_0$ depends only on $(C,[v_C],[\tau|_C])$, and its class
lies in $\mathcal M_b$.
Its image closure $\Delta_b$ is irreducible of dimension $D_b-1$,
hence is a divisor in the $D_b$-dimensional component
$\mathcal M_b$. The image contains a dense open subset of
$\Delta_b$. Its members lie outside the original family by the
uniqueness of the distinguished inclusion and the presence of
vertical fibers in $Z$. The two cited propositions give all the
remaining assertions on this dense locus.
For later use, we note that $T_\partial$ parametrizes exactly the
extensions of Proposition~\ref{prop:torsion-sheaves} with $\mathcal O_C(\Gamma')\simeq
\mathcal O_C(eH)$, that is, with $K\simeq\mathcal O_C$.
\end{proof}

\subsection{A second component through \texorpdfstring{$\Delta_b$}{Delta b}}\label{subsec:second-component}

Throughout this subsection we keep the assumptions of Theorem~\ref{thm:C}.
The boundary sheaves $F_0$ arise from the quotient
$\mathcal V|_C\twoheadrightarrow\mathcal O_C\oplus\mathcal O_\Gamma(eH)$
with $\Gamma=(v_C)_0$. We now allow $\Gamma$ to be an arbitrary reduced
divisor of degree $\ell$ on $C$, not necessarily cut out by a section of
$\mathcal O_C(eH)$.

Let $C\in|H|$ be smooth, let $\Gamma'\subset C$ be a reduced effective
divisor of degree $\ell=2eh$, and put
\[
K=\mathcal O_C(eH-\Gamma'),\qquad
J=\mathcal O_C\oplus\mathcal O_{\Gamma'}(eH),
\]
so that $\deg K=0$ and there is an exact sequence
\[
0\to K\xrightarrow{(0,\iota)}\mathcal O_C\oplus\mathcal O_C(eH)
=\mathcal V|_C\to J\to0,
\]
where $\iota$ is the natural inclusion. With $X_C$, $p$ and $j$
as in the proof of Proposition~\ref{orig:6.12}, set $Z=\PP_C(J)$.

\begin{lemma}\label{lem:torsion-curve}
\begin{enumerate}[label=\textup{(\roman*)}]
\item $Z\subset X_C$ is an effective Cartier divisor with
$\mathcal I_{Z/X_C}\simeq\mathcal O_{X_C}(-\zeta)\otimes p^*K$. It is the union
of the horizontal section defined by
$\mathcal V|_C\twoheadrightarrow\mathcal O_C$ and the fibers
$p^{-1}(q)$, $q\in\Gamma'$, each with multiplicity one; these meet in
ordinary nodes. In particular $Z\subset X$ is a local complete
intersection curve, and $Rp_*\mathcal O_Z\simeq\mathcal O_C$.
\item With $\mathcal L_K$ as in Proposition~\ref{orig:6.12}, one has
\[
\begin{aligned}
\operatorname{Ext}^1_X(\mathcal I_Z(\zeta-bA),\mathcal O_X(bA))
&\simeq H^0(C,\mathcal L_K),\qquad h^0(C,\mathcal L_K)=R_b,\\
\operatorname{Ext}^2_X(\mathcal I_Z(\zeta-bA),\mathcal O_X(bA))&=0.
\end{aligned}
\]
\item There is an isomorphism
$\shExt^1_X(\mathcal I_Z(\zeta-bA),\mathcal O_X(bA))\simeq
p^*\mathcal L_K|_Z$, under which the local class of the extension corresponding to
$\sigma\in H^0(C,\mathcal L_K)$ is $p^*\sigma|_Z$.
\end{enumerate}
\end{lemma}

\begin{proof}
For \textup{(i)}, we use the description in the proof of
Proposition~\ref{orig:6.12}: $Z$ is the zero divisor of the
injective map of line bundles $p^*K\to\mathcal O_{X_C}(\zeta)$ on the
integral surface $X_C$. Near a point of $p^{-1}(q)$, $q\in\Gamma'$, choose
a local parameter $t$ of $C$ at $q$ and a fiber coordinate $x$ vanishing
on the horizontal section; the equation becomes $tx=0$. This gives the
description of the components and the nodes. Since $X_C\subset X$ is
Cartier, $Z\subset X$ is a local complete intersection of codimension
two. Applying $Rp_*$ to
$0\to\mathcal O_{X_C}(-\zeta)\otimes p^*K\to\mathcal O_{X_C}\to\mathcal O_Z\to0$
and using $Rp_*\mathcal O_{X_C}(-\zeta)=0$ gives
$Rp_*\mathcal O_Z\simeq\mathcal O_C$.

Assertion \textup{(ii)} is Proposition~\ref{orig:6.12}\textup{(ii)},
which allows torsion quotients with kernel of degree zero.

For \textup{(iii)}, the local Koszul calculation in the proof of
Proposition~\ref{orig:6.13}\textup{(ii)} applies because $Z\subset X$
is a codimension-two local complete intersection by \textup{(i)}.
It gives
$\shExt^1_X(\mathcal I_Z,\mathcal O_X)\simeq\det N_{Z/X}$. From
$Z\subset X_C\subset X$,
\[
\det N_{Z/X}\simeq N_{Z/X_C}\otimes N_{X_C/X}|_Z
\simeq\mathcal O_Z(\zeta+A)\otimes p^*K^{-1}.
\]
Twisting by $\mathcal O_Z(2bA-\zeta)$ yields $p^*\mathcal L_K|_Z$. Since
$\shHom_X(\mathcal I_Z(\zeta-bA),\mathcal O_X(bA))
=\mathcal O_X(2bA-\zeta)$ has vanishing cohomology in all degrees
(Lemma~\ref{orig:2.8}), the local-to-global spectral sequence identifies
the global extension group with
$H^0(Z,p^*\mathcal L_K|_Z)$. By \textup{(i)} and the projection
formula, this space is $H^0(C,\mathcal L_K)$.
The edge map sends an extension to its
local extension classes; the Koszul identification is independent of
local generators because their change-of-basis determinant is exactly
the transition function of $\det N_{Z/X}$.
\end{proof}

\begin{proposition}\label{prop:torsion-sheaves}
Let $\sigma\in H^0(C,\mathcal L_K)$ be nonzero at every point of $\Gamma'$, with
reduced zero divisor, and let
\[\tag{7.5}\label{eq:torsion-ext}
0\to\mathcal O_X(bA)\to F\to\mathcal I_Z(\zeta-bA)\to0
\]
be the corresponding extension. Such $\sigma$ form a nonempty open
subset of $\PP H^0(C,\mathcal L_K)$. Then:
\begin{enumerate}[label=\textup{(\roman*)}]
\item $F$ is reflexive, locally free along every vertical
component of $Z$, and its singularity scheme is reduced of length
$\deg \mathcal L_K=(4b+2)h$;
\item $F$ has the numerical Chern character $w=\operatorname{ch}(F_0)$;
\item $F$ is slope-stable for every $L\in\mathcal C_b$;
\item $h^0(F)=k_b$ and $h^0(F(-bA))=1$.
\end{enumerate}
\end{proposition}

\begin{proof}
Since $\deg\mathcal L_K=(4b+2)h>2g(C)$, Serre duality gives
$H^1(C,\mathcal L_K(-2q))=0$ for every $q\in C$. The restriction sequence
\[
0\to\mathcal L_K(-2q)\to\mathcal L_K
\to\mathcal L_K\otimes(\mathcal O_C/\mathfrak m_q^2)\to0
\]
therefore induces a surjection on global sections. After choosing a
local parameter at $q$ and a local trivialization of $\mathcal L_K$,
reduction modulo $\mathfrak m_q^2$ records the constant and linear
coefficients of a section. Surjectivity means that these two
coefficients, or equivalently its value and first derivative at $q$,
can be prescribed independently.
In particular, vanishing to order at least two at a fixed $q$ imposes
two independent linear conditions. The incidence locus of pairs
$(q,[\sigma])$ with $\operatorname{ord}_q(\sigma)\ge2$ is consequently
closed of dimension $h^0(C,\mathcal L_K)-2$ in
$C\times\PP H^0(C,\mathcal L_K)$. Since $C$ is projective, its image
in $\PP H^0(C,\mathcal L_K)$ is closed and proper. A general section
therefore has reduced zero divisor. Avoiding also the finitely many
proper hyperplanes $\{\sigma(q)=0\}$, $q\in\Gamma'$, gives the required
nonempty open subset.

For \textup{(i)}, Lemma~\ref{lem:torsion-curve}\textup{(iii)} shows that on each vertical component
the local extension class is the restriction of $p^*\sigma$, which is a
nonzero constant there, because $p^*\mathcal L_K$ is trivial on the fibers of $p$
and $\sigma(q)\ne0$. The local class generates
$\shExt^1$ along the whole vertical component, including the node,
and $F$ is locally free there by Lemma~\ref{orig:4.3}.
That lemma requires only a regular sequence for $\II_Z$, so it
also applies at the nodes of $Z$. Along the horizontal component the local
class vanishes exactly on $(\sigma)_0$. As in Lemma~\ref{orig:4.4}, $F$ is
reflexive, and its singularity scheme is $(\sigma)_0$, reduced of
length $\deg \mathcal L_K$.

For \textup{(ii)}, Proposition~\ref{orig:6.12}\textup{(iii)} gives
the numerical Chern character $w$.

To prove \textup{(iii)}, note that $b\ge e+1$ gives $2b>e$.
Every $L\in\mathcal C_b$ therefore satisfies
\eqref{eq:slope-window} by \eqref{eq:orig:5.1}, and Lemma~\ref{lem:subsheaf-test}(a)
applies to the nonzero torsion-free extension~\eqref{eq:torsion-ext}
with $Z\ne\emptyset$.

Finally, Proposition~\ref{orig:6.12}\textup{(i)} gives $h^0(F)=b^2h+2+h^0(S,\mathcal V(-(b+1)H))$,
and the last term vanishes because $e-b-1<0$. For the second equality,
twist the ideal sequence by $-bA$. We have
$H^0(X,\mathcal O_X(\zeta-(2b+1)A))=H^0(S,\mathcal V(-(2b+1)H))=0$ and
$H^0(C,K(-2bH))=0$ since $\deg K(-2bH)<0$. Hence
$H^0(X,\mathcal I_Z(\zeta-2bA))=0$. The sequence
$0\to\mathcal O_X\to F(-bA)\to\mathcal I_Z(\zeta-2bA)\to0$ then gives
$h^0(F(-bA))=1$.
\end{proof}

\begin{proposition}\label{prop:torsion-family}
The sheaves of Proposition~\ref{prop:torsion-sheaves} form a smooth
irreducible flat family $T'$ of dimension $D_b+h$, whose classifying
morphism $f\colon T'\to M$ is injective on geometric points. Hence
the image $\mathcal F_b^{\mathrm{tor}}:=f(T')$ is irreducible of
dimension $D_b+h$, and its closure in $M$ contains $\Delta_b$.
\end{proposition}

\begin{proof}
Let $B'\subset\operatorname{Sym}^\ell_{U_H}\mathscr C$ be the open
subset of reduced divisors. The relative symmetric product is smooth,
hence open over $U_H$, with irreducible fibers and irreducible base.
Thus $B'$ is
smooth and irreducible of dimension $(h+1)+\ell$.
The universal divisor gives a family of curves $\mathscr Z$ as in
Lemma~\ref{lem:torsion-curve}\textup{(i)}. Its equation is a
nonzerodivisor on each fiber of the universal surface; hence
$\mathscr Z$ is flat and is a relative complete intersection in
$X\times B'$. The embedded curve recovers $C$ and $\Gamma'$ as the
images of its horizontal and vertical components.

Since Lemma~\ref{lem:subsheaf-test}\textup{(a)} gives stability for
every nonzero extension, Lemma~\ref{lem:serre-family} applies.
Let $T'$ be the open subset of the resulting projective extension
bundle defined by the conditions of
Proposition~\ref{prop:torsion-sheaves}. That proposition shows it
is nonempty in every fiber. Thus $T'$ is smooth and irreducible,
carries a flat family, and has dimension
$(h+1)+\ell+R_b-1=D_b+h$.
Its classifying morphism is injective on geometric points, so
$\dim\mathcal F_b^{\mathrm{tor}}=\dim T'$.

Under $(C,[v_C])\mapsto(C,(v_C)_0)$, the base $B_\partial$ of
Corollary~\ref{cor:boundary-divisor} identifies with the locally closed
subset of $B'$ where $\mathcal O_C(\Gamma')\simeq\mathcal O_C(eH)$,
equivalently $K\simeq\mathcal O_C$. For a degree-zero line bundle,
this is equivalent to $h^0(K)\ge1$, a closed condition by upper
semicontinuity. On this subset the extension construction reduces to that of the corollary. The open conditions
on the extension section are also the same: $\Xi=(\sigma)_0$ is
reduced and disjoint from $\Gamma'$. Thus $T_\partial$ identifies
with the restriction of $T'$ to $B_\partial$. Its image is dense
in $\Delta_b$ and lies in $\mathcal F_b^{\mathrm{tor}}$.
Taking closures gives $\Delta_b\subset\overline{\mathcal F_b^{\mathrm{tor}}}$.
\end{proof}

\Needspace{4\baselineskip}
\begin{corollary}\label{cor:second-component}
Let $\mathcal N_b=\overline{\mathcal F_b^{\mathrm{tor}}}$ with its reduced structure.
\begin{enumerate}[label=\textup{(\roman*)}]
\item $\mathcal N_b$ is an irreducible component of $M_{\mathrm{red}}$
of dimension $D_b+h$, distinct from $\mathcal M_b$.
\item $\Delta_b\subset\mathcal M_b\cap\mathcal N_b$,
and $\Delta_b$ is an irreducible component of this intersection.
\item At a general closed point $[F_0]\in\Delta_b$,
\[
\dim\mathcal O_{M,[F_0]}=D_b+h .
\]
\end{enumerate}
\end{corollary}

\begin{proof}
Since $M$ is noetherian, the irreducible closed subset
$\overline{\mathcal F_b^{\mathrm{tor}}}$ lies in some irreducible component $\mathcal N'$ of
$M_{\mathrm{red}}$. Because
$\dim \mathcal N'\ge D_b+h>D_b=\dim\mathcal M_b$, we have
$\mathcal N'\ne\mathcal M_b$. By Proposition~\ref{prop:torsion-family},
$\Delta_b\subset\mathcal M_b\cap \mathcal N'$.

Let $[F_0]\in\Delta_b$ be a general closed point, so that Proposition~\ref{cor:k3-boundary-nilpotents}
applies, and put $R=\mathcal O_{M,[F_0]}$ and
$R_B=\mathcal O_{\mathcal B_b,[F_0]}$. By Proposition~\ref{cor:k3-boundary-nilpotents}(i), $\operatorname{Spec}R$
and $\operatorname{Spec}R_B$ have the same underlying space, so
$\dim R=\dim R_B$ and the two rings have the same minimal primes. By
Proposition~\ref{cor:k3-boundary-nilpotents}(ii), the embedding dimension of $R_B$ is
$\dim T_{[F_0]}\mathcal B_b=\dim\ker\rho=D_b+h+1$. Hence
$\dim R=\dim R_B\le D_b+h+1$. Suppose that
$\dim R=D_b+h+1$. Then $R_B$ is regular, hence a domain, so
$\operatorname{Spec}R$ is irreducible. This contradicts the fact that
the two distinct components $\mathcal M_b$ and $\mathcal N'$ pass
through $[F_0]$. Therefore $\dim R\le D_b+h$. Conversely,
$[F_0]\in\overline{\mathcal F_b^{\mathrm{tor}}}$ and $\overline{\mathcal F_b^{\mathrm{tor}}}$ is irreducible of dimension
$D_b+h$, so $\dim R\ge D_b+h$. This proves (iii).

Since $\mathcal N'$ passes through $[F_0]$, we get
$\dim \mathcal N'\le\dim R=D_b+h=\dim\overline{\mathcal F_b^{\mathrm{tor}}}$. Together with
$\overline{\mathcal F_b^{\mathrm{tor}}}\subset \mathcal N'$ and the irreducibility of both, this forces
$\mathcal N'=\overline{\mathcal F_b^{\mathrm{tor}}}=\mathcal N_b$, proving (i).

For (ii), since neither of the two distinct irreducible components
contains the other, $\mathcal M_b\cap\mathcal N_b$ is a proper
closed subset of $\mathcal M_b$ and has dimension at most
$D_b-1$. Since $\Delta_b$ is irreducible of dimension $D_b-1$
(Corollary~\ref{cor:boundary-divisor}), it is a component of the intersection.
\end{proof}

\begin{remark}
Corollary~\ref{cor:second-component}(iii) sharpens the bound
$\dim R\le D_b+h+1$ of Proposition~\ref{cor:k3-boundary-nilpotents}(ii). In terms of quotients, the
locus of Section~\ref{sec:quot} consists of line-bundle quotients of
$\mathcal V|_C$, of total dimension $2eh+1$. In the full relative Quot
scheme with kernels of degree zero, the torsion quotients
$\mathcal V|_C\twoheadrightarrow J$ form a family of dimension
$2eh+h+1$, and the boundary $\Delta_b$ lies over the intersection of their
closures.
\end{remark}

\Needspace{7\baselineskip}
We now combine the boundary construction, the local Brill--Noether
calculation, and the existence of the second component to prove
Theorem~\ref{thm:C}.

\begin{proof}[Proof of Theorem~\ref{thm:C}]
Proposition~\ref{orig:9.6} gives stability, reflexivity, specialization
to $F_0$, and the cohomology and tangent dimensions.
The component is generically smooth by Corollary~\ref{orig:6.7}.
Proposition~\ref{cor:k3-boundary-nilpotents} identifies the two local
BN ideals and excludes $W^{k_b+1}$ in an ambient neighborhood.
Finally, Corollary~\ref{cor:boundary-divisor} proves the divisor
assertion and extends these conclusions to its general point.
Part (iv) is Corollary~\ref{cor:second-component}, with the description
of the general member given by Propositions~\ref{prop:torsion-sheaves}
and~\ref{prop:torsion-family}.
\end{proof}

\medskip
\noindent\textbf{Acknowledgments.}
The author is deeply grateful to his advisor, Professor Rong Du, for his
insightful advice, valuable discussions, and continuous support during
the writing of this paper.
The author is sponsored by National Natural Science Foundation of China
(Grant No.~12471040) and Science and Technology Commission of Shanghai
Municipality (Grant No.~22DZ2229014).

\phantomsection

\addcontentsline{toc}{section}{References}
\begin{thebibliography}{99}
\small
\setlength{\itemsep}{1pt}\setlength{\parsep}{0pt}
\bibitem{Arrondo} E. Arrondo, A home-made Hartshorne--Serre
correspondence, \textit{Rev. Mat. Complut.} \textbf{20} (2007), 423--443.
\bibitem{SGA6} P. Berthelot, A. Grothendieck and L. Illusie (eds.),
\textit{Th\'eorie des intersections et th\'eor\`eme de Riemann--Roch
(SGA 6)}, Lecture Notes in Mathematics 225, Springer, 1971.
\bibitem{Bridgeland} T. Bridgeland, Equivalences of triangulated categories
and Fourier--Mukai transforms, \textit{Bull. London Math. Soc.}
\textbf{31} (1999), 25--34; corrected version,
\href{https://arxiv.org/abs/math/9809114v2}{arXiv:math/9809114v2}.
\bibitem{BudurWang} N. Budur and B. Wang, Cohomology jump loci of
differential graded Lie algebras, \textit{Compos. Math.}
\textbf{151} (2015), 1499--1528.
\bibitem{RuledSurfaces} L. Costa and I. Mac\'ias Tarr\'io,
Brill--Noether theory of stable vector bundles on ruled surfaces,
\textit{Mediterr. J. Math.} \textbf{21} (2024), Article 118.
\bibitem{RuledThreefolds} L. Costa and I. Mac\'ias Tarr\'io,
Moduli spaces of stable bundles on ruled 3-folds and Brill--Noether
problems, \textit{Ricerche Mat.} \textbf{75} (2026), 2179--2202.
\bibitem{CostaMiroRoig} L. Costa and R. M. Mir\'o-Roig,
Brill--Noether theory for moduli spaces of sheaves on algebraic
varieties, \textit{Forum Math.} \textbf{22} (2010), 411--432.
\bibitem{EJ08} A.-S. Elsenhans and J. Jahnel,
K3 surfaces of Picard rank one and degree two, in \textit{Algorithmic
Number Theory (ANTS VIII)}, Lecture Notes in Computer Science
\textbf{5011}, Springer, Berlin, 2008, 212--225.
\bibitem{Fulton} W. Fulton, \textit{Intersection Theory}, Ergebnisse
der Mathematik und ihrer Grenzgebiete (3), vol. 2, Springer, 1984.
\bibitem{GK} A. Gholampour and M. Kool, Stable reflexive sheaves and
localization, \textit{J. Pure Appl. Algebra} \textbf{221} (2017), 1934--1954.
\bibitem{HartshorneRD} R. Hartshorne, \textit{Residues and Duality},
Lecture Notes in Mathematics 20, Springer, 1966.
\bibitem{HartshorneAG} R. Hartshorne, \textit{Algebraic Geometry},
Graduate Texts in Mathematics 52, Springer, 1977.
\bibitem{Hartshorne} R. Hartshorne, Stable reflexive sheaves,
\textit{Math. Ann.} \textbf{254} (1980), 121--176.
\bibitem{HuybrechtsFM} D. Huybrechts, \textit{Fourier--Mukai Transforms
in Algebraic Geometry}, Oxford Mathematical Monographs,
Clarendon Press, Oxford, 2006.
\bibitem{HuybrechtsK3} D. Huybrechts, \textit{Lectures on K3 Surfaces},
Cambridge Studies in Advanced Mathematics \textbf{158},
Cambridge University Press, Cambridge, 2016.
\bibitem{HL} D. Huybrechts and M. Lehn, \textit{The Geometry
of Moduli Spaces of Sheaves}, 2nd ed., Cambridge Mathematical
Library, Cambridge University Press, Cambridge, 2010.
\bibitem{JMT} M. Jardim, D. Markushevich and A. S. Tikhomirov,
Two infinite series of moduli spaces of rank 2 sheaves on $\PP^3$,
\textit{Ann. Mat. Pura Appl.} (4) \textbf{196} (2017), 1573--1608.
\bibitem{Kleppe} J. O. Kleppe, Moduli spaces of reflexive sheaves
of rank 2, \textit{Canad. J. Math.} \textbf{62} (2010), 1131--1154.
\bibitem{Lazarsfeld} R. Lazarsfeld, Brill--Noether--Petri without
degenerations, \textit{J. Differential Geom.} \textbf{23} (1986), 299--307.
\bibitem{Lieblich} M. Lieblich, Moduli of complexes on a proper
morphism, \textit{J. Algebraic Geom.} \textbf{15} (2006), 175--206.
\bibitem{Mumford} D. Mumford, Further pathologies in algebraic geometry,
\textit{Amer. J. Math.} \textbf{84} (1962), 642--648.
\bibitem{MumfordAV} D. Mumford, \textit{Abelian Varieties},
Tata Institute of Fundamental Research Studies in Mathematics 5,
Oxford University Press, 1970.
\bibitem{TT} R. W. Thomason and T. Trobaugh, Higher algebraic
$K$-theory of schemes and of derived categories, in \textit{The
Grothendieck Festschrift}, Vol.~III, Progress in Mathematics 88,
Birkh\"auser, 1990, 247--435.
\bibitem{VermeireStable} P. Vermeire, Stable reflexive sheaves on
smooth projective 3-folds, \textit{Pacific J. Math.} \textbf{219}
(2005), 391--398.
\bibitem{VermeireModuli} P. Vermeire, Moduli of reflexive sheaves
on smooth projective 3-folds, \textit{J. Pure Appl. Algebra}
\textbf{211} (2007), 622--632.
\end{thebibliography}
\end{document}